\documentclass[journal,twoside,web]{ieeecolor}
\usepackage{generic}
\usepackage{cite}
\usepackage{amsmath,amssymb,amsfonts}
\usepackage{algorithmic}
\usepackage{multirow}
\usepackage{graphicx}
\usepackage{algorithm,algorithmic}
\usepackage{hyperref}
\usepackage{epsfig}
\usepackage{textcomp}
\usepackage{lipsum}

\newtheorem{theorem}{Theorem}[section]
\newtheorem{lemma}[theorem]{Lemma}

\newtheorem{assumption}[theorem]{Assumption}

\newcommand{\dist}{{\rm dist}}
\newcommand{\red}{\textcolor{black}}
\newcommand{\blue}{\textcolor{black}}

\DeclareMathOperator*{\argmax}{arg\,max}

\def\BibTeX{{\rm B\kern-.05em{\sc i\kern-.025em b}\kern-.08em
    T\kern-.1667em\lower.7ex\hbox{E}\kern-.125emX}}
\begin{document}

\title{  Stochastic halfspace approximation method for convex optimization   with nonsmooth  functional constraints}
\author{Nitesh Kumar Singh and Ion Necoara
\thanks{Manuscript received 11 October 2023; revised 6 March 2024 and 29 June 2024; accepted 9 July 2024. Recommended by Associate Editor A. Olshevsky. The research leading to these results has received funding from: UEFISCDI, Romania, PN-III-P4-PCE-2021-0720, project L2O-MOC, nr. 70/2022. (\textit{Corresponding author: Ion Necoara}.)}
\thanks{Nitesh Kumar Singh is with the Automatic Control and Systems Engineering Department, National University of Science and Technology Politehnica Bucharest,   Spl. Independentei 313, 060042 Bucharest, Romania (e-mail: nitesh.nitesh@stud.acs.upb.ro). }
\thanks{Ion Necoara is with the Automatic Control and Systems Engineering Department,  National University of Science and Technology Politehnica Bucharest,   Spl. Independentei 313, 060042 Bucharest, Romania, and with the Gheorghe Mihoc-Caius Iacob Institute of Mathematical Statistics and Applied Mathematics of the  Romanian Academy, 050711 Bucharest, Romania (e-mail: ion.necoara@upb.ro).}
}

\maketitle

\begin{abstract}
In this work, we consider convex optimization problems with smooth objective function and nonsmooth functional constraints. We propose a new stochastic gradient algorithm, called  Stochastic Halfspace Approximation Method (SHAM), to solve this problem, where at each iteration we first take a gradient step for the objective function and then we perform a projection step onto one halfspace approximation of a randomly chosen constraint. We propose various strategies to create this stochastic halfspace approximation and we provide a unified convergence analysis that yields new convergence rates for SHAM algorithm in both optimality and feasibility criteria evaluated at some average point. In particular, we derive convergence rates of order   $\mathcal{O} (1/\sqrt{k})$, when the objective function is only convex, and $\mathcal{O} (1/k)$ when the objective function is strongly convex. The efficiency of SHAM  is illustrated through detailed numerical simulations. 
\end{abstract}

\vspace{-0.15cm}
\section{Introduction}
\label{sec:introduction}
\IEEEPARstart{C}{onvex} optimization, with its wide-ranging applications across various domains such as systems and control \cite{BerBac:22,NedDin:14,NecNed:21},  machine learning\cite{Vap:98},\cite{BhaGra:04}, signal processing \cite{Nec:20}, \cite{Tib:11}, and operations research \cite{RocUry:00}, has evolved to address increasingly complex challenges. Among these challenges, optimization problems with nonsmooth functional constraints stands out as a critical frontier, where traditional approaches face hard limitations when handling a large number of constraints and/or when projections  are computationally demanding. In this context, we address the following convex problem with smooth objective and  nonsmooth functional constraints:
\begin{equation}\label{eq:prob}
    \begin{array}{rl}
	f^* = & \min\limits_{x \red{\in \mathcal{Y}} \in \mathbb{R}^n} \; f(x) \\
	& \text{subject to } \;  h_j(x) \le 0 \;\; \forall j \in [m],
    \end{array}
\end{equation}

\noindent where $[m] = \{ 1,\ldots,m \}$ and the objective function $f:\mathbb{R}^n \rightarrow \mathbb{R}$ and the constraints $h_j:\mathbb{R}^n \rightarrow \mathbb{R}$ are proper, lower-semicontinuous, convex functions. Additionally, $f$ is $L_f$-smooth, while $h_j'$s are general (possibly nondifferentiable) \red{and $\mathcal{Y}$ is a non-empty,  closed, convex and  simple set (by simple we mean that it admits easy projection)}.  Solving problem \eqref{eq:prob} efficiently is pivotal for addressing real-world applications, particularly when $m$ is large, i.e., dealing with a large number of constraints,  and/or when the projection onto some of the subsets $\{ x \in \mathbb{R}^n : h_j (x)\le 0 \}$ are difficult to compute. These issues can be highly computationally demanding, if not impossible. In response to these intricate challenges, the field has witnessed a prevalent adoption of stochastic first order (sub)gradient methods. Notably, the stochastic gradient descent  \cite{Nec:20,GowRic:19} and the stochastic proximal/projection \cite{NecNed:21, Nec:20, WanBer:16} approaches, have emerged as valuable tools in this context (when the projection onto individual sets is easy to compute). Additionally, the realm of problem \eqref{eq:prob} has also seen the introduction of primal-dual stochastic gradient methods, see e.g., \cite{Xu:20}.  
These stochastic optimization algorithms try to strike a balance between computational efficiency, solution quality, and handling many functional constraints.  
Thus the central theme of our scientific inquiry also revolves around harnessing the power of randomness and incorporating stochastic (sub)gradient techniques to tackle problem \eqref{eq:prob}, particularly, those involving many nonsmooth functional constraints. 

\vspace{0.05cm} 

\noindent \textbf{Related work.} In previous studies, e.g.,  \cite{Ned:11}, \cite{NecSin:22}, \cite{SinNec:23},  \cite{YanYue:22} the potential advantages of randomness in enhancing optimization algorithms performance have been already proven. In particular, in \cite{Ned:11} asymptotic convergence results were derived under assumptions analogous to those governing problem \eqref{eq:prob}, while \cite{NecSin:22} offered valuable insights into addressing composite convex optimization problems with nonsmooth convex functional constraints, accompanied by sublinear convergence rates. However, these works either do not provide convergence rates (e.g., \cite{Ned:11} gives only asymptotic convergence for smooth functional constraints case) or require strong assumptions (e.g., compactness of the feasible set  \cite{NecSin:22}). Moreover, both \cite{Ned:11} and \cite{NecSin:22} shared a similar algorithm, and its mini-batch variant was subsequently presented in \cite{SinNec:23}. Furthermore, paper \cite{YanYue:22} considered the problem \eqref{eq:prob} with finite sum objective and developed a variance reduced type method which shares the same convergence behavior as in \cite{Ned:11,NecSin:22}. In this paper, we introduce a novel algorithm, Stochastic Halfspace Approximation Methods (SHAM), tailored to solve problem \eqref{eq:prob}. In each iteration, SHAM executes a gradient step on the objective to yield an intermittent point, followed by projecting this point onto a linear halfspace approximation of a randomly selected constraint with a predefined probability. Importantly, our approach provides a general interpretation for constructing this linear halfspace approximation, applicable to numerous points (in fact, all points on a \blue{line segment}) where we can build this approximation. We offer two distinct examples of this interpretation, accompanied by a unified convergence analysis.


\medskip 

\noindent \textbf{Contributions.} In this paper, 
we extend and synthesize insights from prior research to offer the following contributions:
\\
$(i)$ We introduce a novel Stochastic Halfspace Approximation Method (SHAM) for solving convex optimization problems with smooth objective and nonsmooth functional constraints, which consists of  a gradient step on the objective  followed by a projection onto a linear halfspace approximation of a randomly selected constraint. Our analysis offers various (in fact infinite) examples for creating the stochastic halfspace approximation, accompanied by a unified convergence analysis.\\
$(ii)$ Using new convex tools, we derive sublinear convergence rates of order $\mathcal{O}(1/\sqrt{k})$ for convex objective and of order $\mathcal{O}(1/k)$ when the objective function is strongly convex.\\
$(iii)$ For $L_f$-smooth and $\mu$-strongly convex objective, we propose a switching stepsize strategy, determining when to transit from a constant to a non-increasing stepsize, which depends on the constants $L_f$ and $\mu$.\\
Due to SHAM's novelty, we offer distinct proofs under basic assumptions, setting it apart from the existing works \cite{Ned:11},\cite{NecSin:22}. Our numerical experiments, featuring a quadratic objective and second-order cone (SOC) constraints, consistently validate SHAM's efficiency, affirming its practical applicability. 

\medskip 

\noindent \textbf{Content.} The remainder of this paper is organized as follows: In Section II, we establish essential notations and key assumptions. Section III introduces our algorithm SHAM. In Section IV, we delve into the convergence analysis of SHAM algorithm. Finally, in Section V we present numerical results that validate the practical efficacy of our approach.
\vspace{-0.2cm}
\section{Notations and key assumptions}
\noindent Throughout this paper, we use the following notations. We consider the individual sets $\mathcal{X}_j$ and the feasible set of \eqref{eq:prob}: 
\vspace{-0.2cm}
\[ \mathcal{X}_j \!=\! \{ x\in \mathbb{R}^n  : h_j(x) \le 0\} \;\;\forall j \in [m], \;\; \mathcal{X} \!=\!  \red{\mathcal{Y} \cap (\cap_{j \in [m]} \mathcal{X}_j)}.\]
Further, $f^*$ denotes the optimal value, $\mathcal{X}^*$ the optimal set:
\vspace{-0.2cm}
\[ f^* =  \min_{x\in \mathcal{X}} f(x), \quad \mathcal{X}^* = \{ x \in \mathcal{X} | f(x) = f^* \} \neq \emptyset. \]
Moreover, for any $x \in \mathbb{R}^n$ we denote its projection onto the optimal set $\mathcal{X}^* $ by $\bar{x}$, i.e., $\bar{x} = \Pi_{\mathcal{X}^*}(x)$. For a given scalar $a$, we denote   $(a)_+ = \max (a,0)$. For any $x \in \mathbb{R}^n$, a subgradient is $\Delta_j^x   \in \partial h_j(x)$, where the subdifferential $\partial h_j(x)$ is either a singleton or a nonempty set for any $j \in [m]$. We  denote the linear approximation of the functional constraint $h_j$ at $x$ as:
\vspace{-0.2cm}
 $$l_{h_j}(y; x) = h_j(x) + \langle \Delta_j^x, y-x \rangle.$$
Let us also define, for any $x \in \mathbb{R}^n$ and a subgradient  $\Delta_j^x   \in \partial h_j(x)$, the following convex set:
\vspace{-0.2cm}
\begin{align}\label{eq:half_space}
    \mathcal{L}(h_j;x;\Delta_j^x) \!= \!\! \begin{cases}
       \! \{y \in \mathbb{R}^n\!: l_{h_j}(y; x) \le 0\}, \text{if } \Delta_j^x \!\neq\! {0},\\
        \mathbb{R}^n, \hspace{3.08cm} \text{if } \Delta_j^x \!=\! {0}.
    \end{cases}
\end{align}
 Note that, in the case of $\Delta_j^x \neq {0}$, $\mathcal{L}(h_j;x;\Delta_j^x)$ is a half-space. Also in this paper, we use $\dist(x, \mathcal{X}) = \|x - \Pi_\mathcal{X} (x)\|$.
Throughout the paper, we assume that the set $\mathcal{X}$ is non-empty, 
and the optimization problem \eqref{eq:prob} has a finite optimum. 
Furthermore, for problem  \eqref{eq:prob}, we also assume the following assumptions on the objective function~$f$ and the functional constraints $h_j$'s.

\begin{assumption}
\label{assumption1}
$(i)$ The function $f$ is $\mu$-strongly convex, i.e., there exists $\mu \ge 0$ such that we have the following:
\vspace{-0.3cm}
    \begin{align}\label{eq:strg_conv_f}
        f(y) \!\ge \!f(x) \!+\! \langle \nabla f(x), y -x \rangle + \frac{\mu}{2} \|y - x\|^2 \; \forall x,y \in \red{\mathcal{Y}}.
    \end{align}
$(ii)$ The function $f$ has $L_f$-Lipschitz continuous gradient, i.e., there exists $L_f > 0$ such that for any  $x,y \in \red{\mathcal{Y}}$ it holds:
\vspace{-0.25cm}
    \begin{align}\label{eq:smooth_f}
        f(y) \le f(x) + \langle \nabla f(x), y-x \rangle + \frac{L_f}{2} \|y-x\|^2.
    \end{align}
\end{assumption}
Note that Assumption \ref{assumption1}$(i)$ implies that the function $f$ is convex when $\mu = 0$. Assumption \ref{assumption1}$(ii)$ is a general property of a smooth function (see Lemma $1.2.3$ in \cite{Nes:18}).
\noindent Next, we consider the assumptions for the functional constraints $h_j$.

\begin{assumption}
\label{assumption2}
The convex functional constraints  $h_j'$s satisfy the bounded subgradient condition, i.e., there exists non-negative constant $B_h > 0$ such that for any $x \in \blue{\mathcal{Y}}$ and  $\Delta_j^x \in \partial h_j(x)$, we have:
\vspace{-0.2cm}
\begin{align}
\label{eq:non_smooth_h}
    \|\Delta_j^x\|\le B_h \;\; \forall j \in [m].
\end{align}
\end{assumption}
Finally, our next assumption is a linear regularity type condition for the functional constraints.

\begin{assumption}
\label{assumption3}
    The functional constraints satisfy a regularity condition, i.e., there exists a constant $c>0$ such that:
    \begin{equation}\label{eq:lin_reg}
        \dist(x, \mathcal{X}) \le c \max_{j\in [m]} \left[ (h_j(x))_+ \right]\;\; \forall x \in \blue{\mathcal{Y}}. 
    \end{equation}
\end{assumption}
Note that Assumptions \ref{assumption2} and \ref{assumption3} have been frequently used  in the context of stochastic optimization problems, see e.g.,  \cite{Ned:11},\cite{NecSin:22}.  Particularly, Assumption \ref{assumption3} holds e.g., when the feasible set $\mathcal{X}$ has nonempty interior  or  is polyhedral (see   \cite{LewPan:98}). It also holds for more general sets, e.g., when  the collection of  functional constraints satisfies a strengthened Slater condition, such as the generalized Robinson condition, as detailed in \cite{LewPan:98}[Corollary 3]. Next, we provide a lemma that has an important role in proving the convergence results later. The proof follows similarly as in \cite{YanYue:22}(Proposition $1(ii)$). However, instead of using the compactness of a set $\mathcal{Y}$ as in \cite{YanYue:22}, we use Assumption \ref{assumption2}. 
\begin{lemma}\label{lemma_lin_reg}
    Let us assume that $h_j$'s are convex functions and that the feasible set $\mathcal{X}$ is non-empty. Additionally, we assume that Assumptions \ref{assumption2} and \ref{assumption3} hold. Then, the following relation is valid for any $x \in \blue{\mathcal{Y}}$:
    \vspace{-0.2cm}
    \begin{align}\label{eq:lin_reg2}
        & \dist(x, \mathcal{X}) \le c B_h \max_{j \in [m]} \min_{\Delta_j^x \in \partial h_j(x)} \dist(x, \mathcal{L} (h_j;x;\Delta_j^x)). 
    \end{align}
\end{lemma}

\begin{proof}
See Appendix. 
\end{proof}

\vspace{-0.3cm}
\section{Algorithm}
\noindent In this section, we introduce a novel Stochastic Halfspace Approximation Method (SHAM), in order to solve problem \eqref{eq:prob}. For any given iteration  $k$, we choose a random index $j_k \in [m]$ with given probability $\mathbb{P}$, and update as:
\begin{algorithm}[H]
\caption{(SHAM)}\label{alg:alg1}
\begin{algorithmic}
\STATE $\text{Choose} \; x_0 \in \red{\mathcal{Y}} \; \text{stepsize sequence} \; \{\alpha_k\}_{k\geq 0}$  and $\beta > 0$. 
\STATE For {$k \geq 0$  update:}
\vspace{-0.3cm}
\begin{align}
    & \blue{u_k = } x_{k} - \alpha_k \nabla f(x_{k}) \label{eq:alg2step1}\\[-3pt]
    & v_k = \blue{\Pi_\mathcal{Y} (u_k)} \label{eq:alg2step}\\[-3pt]
    &  \text{Sample} \; j_k \sim    \mathbb{P}, \text{choose } \tilde{x}_k,
     \Delta_{j_k}^{\tilde{x}_k} \in \partial h_{j_k} (\tilde{x}_k)\; \text{and update:} \vspace{-0.3cm} \nonumber  \\[-3pt]
    & \red{z_k = (1 - \beta) v_k + \beta \Pi_{\mathcal{L}_{j_k}} (v_k)} \vspace{-0.3cm} \label{eq:alg2step2}\\[-3pt]
    & \red{x_{k+1} = \Pi_\mathcal{Y} (z_k).} \label{eq:alg2step3}
\end{align}
\end{algorithmic}
\end{algorithm}

\vspace{-0.2cm}

\noindent Here  $\mathcal{L}_{j_k} = \mathcal{L}(h_{j_k}; \tilde{x}_k; \Delta_{j_k}^{\tilde{x}_k})$ (see \eqref{eq:half_space}), and the choice of $\tilde{x}_k\in \blue{\mathcal{Y}}$ is user dependent. \textit{Note that the algorithms in \cite{Ned:11,NecSin:22} consider only the case  $\tilde{x}_k = v_k$, whereas SHAM introduces a novel scheme in comparison to these existing works, allowing us the flexibility to create linear halfspace approximations at any point $\tilde{x}_k$ selected on the line \blue{segment} joining the two points $x_k$ and $v_k$. More precisely, we conduct below the convergence analysis for a general choice of $\tilde{x}_k = \gamma v_k + (1 - \gamma) x_k$, for  $\gamma \in \blue{[0,\; 1]}$, that particularly takes the form $\tilde{x}_k = x_k$ when $\gamma = 0$ or $ \tilde{x}_k = v_k$ when $\gamma = 1$, for generating the linear halfspace approximations in step \eqref{eq:alg2step2}. 
Later in Section V, one can also notice the benefits of these choices in practice.} Consequently, SHAM comprises two key steps: firstly, step \eqref{eq:alg2step1} is the gradient step on the objective function seeking to reduce $f$ 
that together with the projection step \eqref{eq:alg2step} gives us an intermittent point $v_k$.
Next, in step \red{\eqref{eq:alg2step2}} we have a  (sub)gradient projection step striving to minimize the feasibility violation of the chosen random constraint $j_k \in [m]$. In this step, if $\Delta_{j_k}^{\tilde{x}_k} \neq 0$, we project  $v_k$ onto the linear halfspace $l_{h_{j_k}}(y; \tilde{x}_k) \le 0$. Thus, \eqref{eq:alg2step2}  has the following  expression:
\vspace{-0.3cm}
\begin{align}\label{eq:algstep}
    \red{z_k} = v_k - \beta \frac{(l_{h_{j_k}}(v_k; \tilde{x}_k))_+}{\|\Delta_{j_k}^{\tilde{x}_k}\|^2} \Delta_{j_k}^{\tilde{x}_k},
\end{align} 
\vspace{-0.1cm}

\noindent where $\beta > 0$. While if $\Delta_{j_k}^{\tilde{x}_k} = 0$, \eqref{eq:alg2step2} yields:
\vspace{-0.2cm}
\begin{align}\label{eq:algstep2}
    \red{z_k} = v_k.
\end{align}

\vspace{-0.1cm}

\noindent Note that in both expressions \eqref{eq:algstep} and \eqref{eq:algstep2} we use the formula of projection onto a halfspace.  Moreover, when $\tilde{x}_k = v_k$ the step \eqref{eq:algstep} uses a Polyak  type stepsize of the form $\beta (h_{j_k}(v_k))_+/\|\Delta_{j_k}^{\tilde{x}_k}\|^2$ with a parameter $\beta > 0$, see also  \cite{Pol:69}, \cite{Pol:01}.
In our analysis, we also need an assumption concerning the distribution of constraint index $j_k$.
\begin{assumption}
    \label{assumption4}
    There exists a constant $\rho \in (0,1]$ such that
    \vspace{-0.3cm}
    \[ \inf_{k \ge 0} \mathbb{P}(j_k = j| \mathcal{F}_{[k]}) \ge \frac{\rho}{m} \quad  \text{a.s.} \quad   \forall j \in [m], \]
\vspace{-0.5cm} 

\noindent where $\mathcal{F}_{[k]} = \{ x_0, j_0, ..., j_{k-1} \}$ for $k\ge 1$ and $\mathcal{F}_0 = \{x_0\}$.
\end{assumption}
This assumption is also standard in the literature, see \cite{WanBer:16}, \cite{YanYue:22}.
In our analysis below, we  consider boundedness of the  norm of the gradients of $f$ along the iterates of SHAM, i.e.:
\vspace{-0.1cm}
\begin{align} \label{eq:bound_xk}
    \|\nabla f(x_k)\| \le B_f \quad \forall k\ge 0.
\end{align} 

\vspace{-0.2cm}

\noindent \red{ Note that since the gradient $\nabla f$ is continuous and the update  step \eqref{eq:alg2step3} ensures $x_k \in \mathcal{Y}$,   condition \eqref{eq:bound_xk} always holds provided that  the set $\mathcal{Y}$ is  bounded,   (see also \cite{BooLan:23}). However, boundedness of $\mathcal{Y}$ is a sufficient condition for \eqref{eq:bound_xk} (there may be problems having  unbounded $\mathcal{Y}$ for which \eqref{eq:bound_xk} still holds).}

\vspace{-0.2cm}
\section{Convergence analysis of SHAM} 

\noindent In this section, we derive the convergence properties of the algorithm SHAM. Before providing the main results, we give some useful inequalities which we use throughout our analysis. We start with some  basic properties of the projection onto a closed convex set $\mathcal{Y} \subseteq \mathbb{R}^n$ (see e.g., \cite{NecSin:22}):
\vspace{-0.2cm}
\begin{align}
    & \|p - \Pi_\mathcal{Y} (p)\| \le \|p - q\| \quad \forall p \in \mathbb{R}^n, q \in \mathcal{Y},\label{eq:proj_property}\\[-3pt]
    & \red{\|q - \Pi_\mathcal{Y} (p)\| \le \|p - q\| \quad \forall p \in \mathbb{R}^n, q \in \mathcal{Y}.}\label{eq:proj_property2}
\end{align}
\vspace{-0.5cm}

\noindent Further, we also have:
\vspace{-0.2cm}
\begin{align}
    \|p+q\|^2 \le 2 \|p\|^2 + 2 \|q \|^2 \quad \forall p, q \in \mathbb{R}^n, \label{eq:ineq2} 
\end{align}
\vspace{-0.5cm}

\noindent and lastly, for any scalar $a$, the following is true:
\vspace{-0.2cm}
\begin{align}\label{eq:ineq3}
    (a)_+ a = (a)_+^2.
\end{align}
\vspace{-0.5cm}

\noindent Now, we provide some lemmas to use later for the convergence results. In the first lemma, we prove a kind of descent property of $f$ using the smoothness condition \eqref{eq:smooth_f}.


\begin{lemma}
\label{lemma_smooth_f}
	Let Assumption \ref{assumption1}$(ii)$ hold. Then, the following relation  is true for the sequences generated by SHAM:
 
 \vspace{-0.6cm}

 \begin{align*}
		& 2\langle \blue{u_k} - x_{k+1}, x_{k+1} - x_k \rangle  \\[-2pt]
        & \le 2 \alpha_k ( f(x_k) - f(x_{k+1})) - (2 - \alpha_k L_f)\|x_{k+1} - x_k\|^2.
	\end{align*}
\end{lemma}

\begin{proof}
See Appendix.
\end{proof}

\noindent Next, using the convexity and other key assumptions on $h_j$'s we develop a relation between the sequence $\tilde{x}_k$ and $v_k$. 

\begin{lemma}\label{lemma_non_smooth_h2}
    Assume that $h_j$'s are convex functions and the feasible set $\mathcal{X}$ is non-empty. Additionally,  Assumptions \ref{assumption2}, \ref{assumption3}, \ref{assumption4} hold. Then, we have the following relation true:
    \vspace{-0.2cm}
    \begin{align*}
        & \dist^2 (\tilde{x}_k, \mathcal{X}) \\
        & \le \frac{2 m c^2 B_h^2}{\rho}  (\mathbb{E}\left[\|v_k - \Pi_{\mathcal{L}_{j_k}} (v_k)\|^2 |\mathcal{F}_{[k]}\right] + \|v_k - \tilde{x}_k\|^2). 
    \end{align*}
\end{lemma}
\begin{proof}
See Appendix.
\end{proof}

\noindent In the following lemma we prove a distance relation for the points $x_{k+1}$ and $v_k$ from the feasible set $\mathcal{X}$.

\begin{lemma}\label{lemma_dist}
Assume that the function $h_j$'s are convex. Then, for $\beta \in (0,2)$, the following relation is true for the iterates generated by SHAM:
\vspace{-0.1cm}
    \[ \dist^2 (x_{k+1}, \mathcal{X}) \le \dist^2 (v_k, \mathcal{X}). \]
\end{lemma}
\begin{proof}
See Appendix.
\end{proof}
\noindent Now,  we prove a relation between  points $x_{k+1}$, $\tilde{x}_k$, and $v_k$.

\begin{lemma}\label{lemma_non_smooth_h}
    Under the assumption of Lemma \ref{lemma_non_smooth_h2} with $\beta \in (0,2)$, the following relation holds:
    \vspace{-0.1cm}
    \begin{align*}
        & \frac{\rho}{4 m c^2 B_h^2}\mathbb{E}[ \dist^2 (x_{k+1}, \mathcal{X})|\mathcal{F}_{[k]}] \\
        & - \!\left(\! 1 \!+\! \frac{\rho}{2 m c^2 B_h^2}\!\right) \! \|v_k - \tilde{x}_k\|^2 \le  \mathbb{E} \! \left[\|v_k - \Pi_{\mathcal{L}_{j_k}} (v_k)\|^2 |\mathcal{F}_{[k]}\right].
    \end{align*}
\end{lemma}

\begin{proof}
See Appendix.
\end{proof}

\noindent Now, we are ready to provide the main recurrence for SHAM. 

 \medskip
 
\begin{theorem}
\label{lemma_common}
    Let us assume that the functions $h_j$'s are convex, the condition \eqref{eq:bound_xk} holds, and the function $f$  satisfies  Assumption \ref{assumption1}$(i)$. Additionally,   Assumptions \ref{assumption1}$(ii)$, \ref{assumption2}, \ref{assumption3} and \ref{assumption4} hold. Then, choosing $\alpha_k \in \left(0, \frac{1}{L_f}\right]$, $\beta \in \left(0, 1 \right)$ and $\tilde{x}_k = \gamma v_k + (1 - \gamma) x_k$ with  \blue{$\gamma \in \blue{[0,\; 1]}$}, the iterates generated by SHAM satisfies the following recurrence: 
   \vspace{-0.2cm}
    \begin{align} \label{eq:common}
        &\mathbb{E}[\|x_{k+1} - \bar{x}_{k+1}\|^2 ] \\
        & \le (1 - \mu \alpha_k)\mathbb{E}[\|x_k - \bar{x}_k\|^2] - 2 \alpha_k \mathbb{E}[ (f(x_{k+1}) - f(\bar{x}_k)) ]\nonumber\\
        & \quad - \frac{\rho \beta (1 - \beta)}{\blue{4} m c^2 B_h^2}\mathbb{E}[ \dist^2 (x_{k+1}, \mathcal{X})] + \alpha_k^2 B^2, \nonumber
    \end{align}
    where $B^2 \!=\! B_f^2\left(  \blue{\frac{1}{1 - \beta} + \beta (1 - \beta)}(1 - \gamma)^2\left( 1 + \frac{\rho}{2 m c^2 B_h^2}\right)\right).$
\end{theorem}

\begin{proof}  Recall that $\bar{x}_k$ is the projection of $x_k$ onto the optimal set $\mathcal{X}^*$, \red{i.e., $\bar{x}_k\in \mathcal{X}^* \subseteq \mathcal{Y}$}. Thus, we have:
\vspace{-0.15cm}
\begin{align*}
	& \|x_{k+1} \!-\! \bar{x}_{k+1}\|^2 \red{\!\overset{\eqref{eq:proj_property}, \eqref{eq:proj_property2}}{\le}\! \|z_k \!-\! \bar{x}_k\|^2 }\!=\! \|\red{z_k} \!-\!v_k \!+\! v_k\! -\! \bar{x}_k\|^2\nonumber\\
	& = \|v_k - \bar{x}_k\|^2 + 2 \langle \red{z_k} -v_k, v_k - \bar{x}_k \rangle + \| \red{z_k} -v_k \|^2\nonumber \\
    & \overset{\eqref{eq:proj_property2}}{\blue{\le}} \|\blue{u_k} - \bar{x}_k\|^2 + 2 \langle \red{z_k} -v_k, v_k - \bar{x}_k \rangle + \| \red{z_k} -v_k \|^2\nonumber \\
    & \overset{\eqref{eq:alg2step1}}{=} \|x_k \!-\! \bar{x}_k\|^2 + 2 \alpha_k \langle \nabla f(x_k), \bar{x}_k\!-\! x_k \rangle \!+\!\|\blue{u_k}- x_k\|^2 \nonumber\\
	& \quad   + 2 \langle \red{z_k} -v_k, v_k - \bar{x}_k \rangle +\| \red{z_k} -v_k \|^2 \nonumber\\
    &  \overset{\eqref{eq:strg_conv_f}}{\le} (1 - \mu \alpha_k)\|x_k - \bar{x}_k\|^2 - 2 \alpha_k ( f(x_k) - f(\bar{x}_k))  \nonumber\\
    & \quad  + \|\blue{u_k}- x_k\|^2  + 2 \langle \red{z_k} -v_k, v_k - \bar{x}_k \rangle +\| \red{z_k} -v_k \|^2 \nonumber\\
    & \quad \!+\! 2\langle \blue{u_k} - x_{k+1}, x_{k+1} - x_k \rangle \!-\! 2\langle \blue{u_k} - x_{k+1}, x_{k+1} - x_k \rangle\nonumber\\
    & \le (1 - \mu \alpha_k)\|x_k - \bar{x}_k\|^2 - 2 \alpha_k ( f(x_{k+1}) - f(\bar{x}_k))\nonumber\\
	& \quad  + \|\blue{u_k}- x_k\|^2 + 2 \langle \red{z_k} -v_k, v_k - \bar{x}_k \rangle +\| \red{z_k} -v_k \|^2\nonumber\\
    & \quad - (2\!-\! \alpha_k L_f)\|x_{k+1} \!\!-\! x_k \|^2 \! -\! 2\langle \red{ z_k} - x_{k+1}, x_{k+1} - x_k \rangle  \\
    & \quad -\! 2\langle \blue{u_k - v_k}, x_{k+1} \blue{-  v_k  + v_k} - x_k \rangle \!-\! 2\langle \blue{ v_k} \red{- z_k}, x_{k+1} - x_k \rangle \nonumber\\
    & \le (1 - \mu \alpha_k)\|x_k - \bar{x}_k\|^2 - 2 \alpha_k ( f(x_{k+1}) - f(\bar{x}_k))\nonumber\\
	& \quad  + \|\blue{u_k}- x_k\|^2 + 2 \langle \red{z_k} -v_k, v_k - \bar{x}_k \rangle + \blue{2} \| \red{z_k} -v_k \|^2\nonumber\\
    & \quad - (2- \alpha_k L_f)\|x_{k+1} - x_k \|^2 \blue{+ \frac{1}{\eta} \| u_k - v_k \|^2} \\
    & \quad \blue{+ \eta \;\| x_{k+1} - v_k \|^2 + \|x_{k+1} - x_k\|^2 }\nonumber\\
    & \overset{\eqref{eq:proj_property}, \eqref{eq:proj_property2}}{\le} (1 - \mu \alpha_k)\|x_k - \bar{x}_k\|^2 - 2 \alpha_k ( f(x_{k+1}) - f(\bar{x}_k)) \\
	& \quad + \blue{\left(1 + \frac{1}{\eta}\right)}\|\blue{u_k}- x_k\|^2  + 2 \langle \red{z_k} -v_k, v_k - \bar{x}_k \rangle \nonumber\\
    & \quad + \blue{\left( 2 + \eta \right)}\| \red{z_k} - v_k \|^2 - \left(1 - \alpha_k L_f \right)\|x_{k+1} - x_k \|^2, \nonumber
\end{align*}


\noindent where the third inequality follows from Assumption \ref{assumption1}$(i)$ (note that if $f$ is only a convex function then this inequality is still true with $\mu = 0$), in the fourth inequality we use Lemma \ref{lemma_smooth_f}, and in fifth inequality we use \red{ the optimality conditions for steps \eqref{eq:alg2step} and \eqref{eq:alg2step3}, i.e., \blue{$\langle u_k \!-\! v_k, v_k \!-\! x_k \rangle \ge 0$} and $\langle z_k \!-\! x_{k+1}, x_{k+1} \!-\! x_k \rangle \ge 0$}, respectively, and \blue{$2\langle p, q\rangle \le \frac{1}{\eta}\|p\|^2 + \eta\|q\|^2 \; \forall p,q \in \mathbb{R}^n$ and  $\forall  \eta > 0$ (with $\eta = 1$ in last term)}. Noticing that $\alpha_k \in \left(0, \frac{1}{L_f}\right]$, we further get:
\begin{align}
    \label{eq:middle1}
    & \|x_{k+1} - \bar{x}_{k+1}\|^2 \le (1 - \mu \alpha_k)\|x_k - \bar{x}_k\|^2\\
    & \quad - 2 \alpha_k ( f(x_{k+1}) - f(\bar{x}_k)) + \blue{\left(1 + \frac{1}{\eta} \right)}\|\blue{u_k}- x_k\|^2 \nonumber\\
    & \quad + 2 \langle \red{z_k} -v_k, v_k - \bar{x}_k \rangle \!+\! \blue{\left( 2 + \eta \right)} \| \red{z_k} -v_k \|^2. \nonumber
\end{align}


\noindent For a given $\Delta_{j_k}^{\tilde{x}_k} \in \partial h_{j_k}( \tilde{x}_k)$, we first consider the case when $\Delta_{j_k}^{\tilde{x}_k} = 0 \overset{\eqref{eq:algstep2}}{\implies} \red{z_k} = v_k$,  and noticing that by the choice of $\mathcal{L}_{j_k}$ from \eqref{eq:half_space}, we have $\Pi_{\mathcal{L}_{j_k}} (v_k) = v_k$, thus  \eqref{eq:middle1} reduces to:


\begin{align}\label{eq:middle2}
    & \|x_{k+1} - \bar{x}_{k+1}\|^2 \le (1 - \mu \alpha_k)\|x_k - \bar{x}_k\|^2  \nonumber\\
    &  \quad - 2 \alpha_k ( f(x_{k+1}) - f(\bar{x}_k)) + \blue{\left(1 + \frac{1}{\eta}\right)}\|\blue{u_k}- x_k\|^2 \nonumber\\
    & \quad - \beta (2 - \blue{\left( 2 + \eta \right)} \beta) \| v_k - \Pi_{\mathcal{L}_{j_k}} (v_k)\|^2. 
\end{align}


\noindent Next, consider the case when $\Delta_{j_k}^{\tilde{x}_k} \neq 0$, from the definition of \red{$z_k$} given by \eqref{eq:algstep}, relation \eqref{eq:middle1} will be:
\vspace{-05cm}
\begin{align*}
    & \|x_{k+1} - \bar{x}_{k+1}\|^2 \le (1 - \mu \alpha_k)\|x_k \!-\! \bar{x}_k\|^2 \!+\! \blue{\left( 2+ \eta \right)} \| \red{z_k} \!-\!v_k \|^2 \nonumber\\[-2pt]
    & \quad - 2 \alpha_k ( f(x_{k+1}) - f(\bar{x}_k)) +\blue{\left(1 + \frac{1}{\eta}\right)}\|\blue{u_k}- x_k\|^2 \nonumber\\[-2pt]
    & \quad + 2\beta \frac{ (l_{h_{j_k}}(v_k; \tilde{x}_k))_+}{\|\Delta_{j_k}^{\tilde{x}_k}\|^2} \langle \Delta_{j_k}^{\tilde{x}_k}, \bar{x}_k -\tilde{x}_k + \tilde{x}_k - v_k \rangle \nonumber\\[-2pt]
    & \le (1 - \mu \alpha_k)\|x_k - \bar{x}_k\|^2 - 2 \alpha_k ( f(x_{k+1}) - f(\bar{x}_k))\nonumber \\[-2pt]
    & \quad +\blue{\left(1 + \frac{1}{\eta}\right)}\|\blue{u_k}- x_k\|^2 + \blue{\left( 2 + \eta \right)} \| \red{z_k} - v_k \|^2 \nonumber\\[-2pt]
    & \quad + 2\beta \frac{ (l_{h_{j_k}}(v_k; \tilde{x}_k))_+}{\|\Delta_{j_k}^{\tilde{x}_k}\|^2} ( h_{j_k}(\bar{x}_k) - h_{j_k}(\tilde{x}_k))\nonumber\\[-2pt]
    & \quad - 2\beta \frac{ (l_{h_{j_k}}(v_k; \tilde{x}_k))_+}{\|\Delta_{j_k}^{\tilde{x}_k}\|^2} \langle \Delta_{j_k}^{\tilde{x}_k}, v_k - \tilde{x}_k \rangle\nonumber\\[-2pt]
    & \le (1 - \mu \alpha_k)\|x_k - \bar{x}_k\|^2 - 2 \alpha_k ( f(x_{k+1}) - f(\bar{x}_k)) \nonumber\\[-2pt]
    & \quad +\blue{\left(1 + \frac{1}{\eta}\right)}\|\blue{u_k}- x_k\|^2 + \blue{\left( 2 + \eta \right)} \beta^2 \frac{ (l_{h_{j_k}}(v_k; \tilde{x}_k))^2_+}{\|\Delta_{j_k}^{\tilde{x}_k}\|^2} \nonumber\\[-2pt]
    & \quad - 2\beta \frac{ (l_{h_{j_k}}(v_k; \tilde{x}_k))_+}{\|\Delta_{j_k}^{\tilde{x}_k}\|^2} l_{h_{j_k}}(v_k; \tilde{x}_k)\nonumber,
\end{align*}


\noindent where the second inequality uses the convexity of $h_j$, and the third inequality uses the fact that the constraints are feasible at point $\bar{x}_k$, i.e., $h(\bar{x}_k)\le 0$. Now, using \eqref{eq:ineq3} and from the expression of \red{$z_k$ in \eqref{eq:alg2step2}} and \eqref{eq:algstep}, we get:
\vspace{-0.2cm}
\begin{align}\label{eq:middle3}
    & \|x_{k+1} - \bar{x}_{k+1}\|^2 \nonumber\\[-2pt]
    & \overset{\eqref{eq:ineq3}}{\le} (1 - \mu \alpha_k)\|x_k - \bar{x}_k\|^2 - 2 \alpha_k ( f(x_{k+1}) - f(\bar{x}_k))\nonumber \\[-2pt]
    & \quad \!+\!\blue{\left(1 \!+\! \frac{1}{\eta}\right)}\|\blue{u_k}\!-\! x_k\|^2 - \beta (2 - \blue{\left( 2 + \eta \right)} \beta) \frac{ (l_{h_{j_k}}(v_k; \tilde{x}_k))^2_+}{\|\Delta_{j_k}^{\tilde{x}_k}\|^2}\nonumber \\[-2pt]
    & = (1 - \mu \alpha_k)\|x_k - \bar{x}_k\|^2 - 2 \alpha_k ( f(x_{k+1}) - f(\bar{x}_k)) \\[-2pt]
    & \quad \!+\!\blue{\left(1 \!+\! \frac{1}{\eta}\right)}\|\blue{u_k}\!-\! x_k\|^2 \!-\! \beta (2 \!-\! \blue{\left( 2 + \eta \right)} \beta) \| v_k \!-\! \Pi_{\mathcal{L}_{j_k}} (v_k)\|^2.\nonumber
\end{align}

\noindent For any given $\Delta_{j_k}^{\tilde{x}_k}$, from relations \eqref{eq:middle2} and \eqref{eq:middle3}, we have:
\vspace{-0.1cm}
\begin{align*}
    & \|x_{k+1} - \bar{x}_{k+1}\|^2\\[-2pt]
    & \le (1 - \mu \alpha_k)\|x_k - \bar{x}_k\|^2 - 2 \alpha_k ( f(x_{k+1}) - f(\bar{x}_k))\nonumber \\[-2pt]
    & \quad \!+\!\blue{\left(1 \!+\! \frac{1}{\eta}\right)}\|\blue{u_k}\!-\! x_k\|^2 \!-\! \beta (2 \!-\! \blue{\left( 2 + \eta \right)} \beta) \| v_k \!-\! \Pi_{\mathcal{L}_{j_k}} (v_k)\|^2.
\end{align*}
\vspace{-0.4cm}

\noindent \blue{By setting $\eta = \frac{1-\beta}{\beta} >0$}, we get:
\vspace{-0.2cm}
\begin{align*}
    & \|x_{k+1} - \bar{x}_{k+1}\|^2\\[-2pt]
    & \le (1 - \mu \alpha_k)\|x_k - \bar{x}_k\|^2 - 2 \alpha_k ( f(x_{k+1}) - f(\bar{x}_k))\nonumber \\[-2pt]
    & \quad \!+\!\blue{\left(\frac{1}{1 - \beta} \right)}\|\blue{u_k}\!-\! x_k\|^2 \!-\! \blue{\beta (1 - \beta)} \| v_k - \Pi_{\mathcal{L}_{j_k}} (v_k)\|^2.
\end{align*}
\vspace{-0.4cm}

\noindent After taking expectation w.r.t. $j_k$ conditioned on $\mathcal{F}_{[k]}$ and using Lemma \ref{lemma_non_smooth_h}, we get:
\vspace{-0.3cm}
\begin{align}\label{eq:middle4}
    & \mathbb{E}[\|x_{k+1} - \bar{x}_{k+1}\|^2 | \mathcal{F}_{[k]} ] \le (1 - \mu \alpha_k)\|x_k - \bar{x}_k\|^2 \nonumber\\[-2pt]
    & \quad - 2 \alpha_k \mathbb{E}[ (f(x_{k+1}) - f(\bar{x}_k))| \mathcal{F}_{[k]} ] \!+\!\blue{\left(\frac{1}{1 - \beta} \right)} \|\blue{u_k} - x_k\|^2 \nonumber \\[-2pt]
    & \quad  - \frac{\rho\beta (1 - \beta)}{\blue{4} m c^2 B_h^2}\mathbb{E}[ \dist^2 (x_{k+1}, \mathcal{X})|\mathcal{F}_{[k]}] \nonumber\\[-2pt]
    & \quad +  \blue{\beta (1 - \beta)}\left( 1 + \frac{\rho}{2 m c^2 B_h^2}\right) \|v_k - \tilde{x}_k\|^2.
\end{align}
Now, using the expression of $\tilde{x}_k$, i.e., $\tilde{x}_k = \gamma v_k + (1 - \gamma)x_k$, for  $\gamma \in \blue{[0, 1]}$, and \eqref{eq:alg2step1}, \eqref{eq:proj_property2} in \eqref{eq:middle4}, we  finally get:
\vspace{-0.3cm}
\begin{align*}
    & \mathbb{E}[\|x_{k+1} - \bar{x}_{k+1}\|^2 | \mathcal{F}_{[k]} ]\nonumber\\[-2pt]
    &  \le (1 - \mu \alpha_k)\|x_k - \bar{x}_k\|^2 - 2 \alpha_k \mathbb{E}[ (f(x_{k+1}) - f(\bar{x}_k))| \mathcal{F}_{[k]} ]\nonumber \\[-2pt]
    & \quad  - \frac{\rho \beta (1 - \beta)}{\blue{4} m c^2 B_h^2}\mathbb{E}[ \dist^2 (x_{k+1}, \mathcal{X})|\mathcal{F}_{[k]}] \nonumber\\[-2pt]
    & \;\; \!+\! B_f^2\left( \blue{\frac{1}{1 - \beta}  + \beta (1 - \beta)}(1 - \gamma)^2\left( 1 + \frac{\rho}{2 m c^2 B_h^2}\right)\right)\alpha_k^2 .
\end{align*}
After taking full expectation we get the required result.
\end{proof}

\noindent It is worth noting that the relationship described in \eqref{eq:common} is also applicable to convex objective functions (i.e., $\mu = 0$). \textit{Moreover, from Lemmas 4.1-4.4 and \eqref{eq:bound_xk} it becomes evident that the derivation of the main recurrence \eqref{eq:common} employs novel and distinct tools compared to those used in \cite{Ned:11, NecSin:22}. Consequently, our convergence analysis introduces new insights and offers different proofs.} With this foundation, we are poised to present the convergence rates of SHAM, which we delineate in the subsequent sections, taking into account the characteristics of the objective function $f$.
\setlength{\belowdisplayskip}{-0.03cm} 
\setlength{\belowdisplayshortskip}{-0.03cm}
\vspace{-0.2cm}

\subsection{Convergence rates  for convex objective function}
\label{Sec_con_smooth_nonsmooth}

\noindent In this section, we consider the scenario where the objective function $f$ in problem \eqref{eq:prob} is only smooth and convex. Under the other assumptions defined in previous sections (except Assumption \ref{assumption1}$(i)$, where $\mu=0$), we provide the convergence rates of Algorithm \ref{alg:alg1}. Before presenting the main results, let us define the following  average sequences:
\vspace{-0.3cm}
\[ \hat{x}_{k} = \frac{1}{S_k}  \sum_{t=0}^{k-1} \alpha_t x_{t+1}, \text{ and } \hat{x}^*_{k} =  \frac{1}{S_k}  \sum_{t=0}^{k-1} \alpha_t\bar{x}_{t},  \]
where $S_k = \sum_{t=0}^{k-1}\alpha_t$.
The following theorem outlines the convergence rates of Algorithm \ref{alg:alg1}.

\begin{theorem}\label{Th:convex_case}
\blue{Under the conditions of Theorem \ref{lemma_common} with convex $f$, i.e., $\mu = 0$,} we have the following bounds  for the average sequence $\hat{x}_{k}$ in terms of optimality and feasibility violation for problem \eqref{eq:prob}:
    \vspace{-0.3cm}
	\begin{align*}
		& \mathbb{E}[(f(\hat{x}_{k}) - f(\hat{x}^*_{k}))] \le \frac{\|x_0 - \bar{x}_0\|^2}{2 S_k} + B^2 \frac{\sum_{t=1}^{k} \alpha_t^2}{2 S_k},\\[-2pt]
		& \mathbb{E}[\dist^2 (\hat{x}_{k}, \mathcal{X})] \\[-2pt]
        & \le \frac{\blue{4}\alpha_0 m c^2 B_h^2}{ \rho \beta (1 - \beta)} \left( \frac{\|x_0 - \bar{x}_0\|^2}{S_k} + B^2 \frac{\sum_{t=1}^{k} \alpha_t^2}{  S_k} \right).
	\end{align*}
\end{theorem}

\begin{proof}
    For convex case, the inequality \eqref{eq:common} holds with $\mu =0$, thus we have:
    \vspace{-0.2cm}
    \begin{align*}
        &\mathbb{E}[\|x_{k+1} - \bar{x}_{k+1}\|^2 ] \\[-1pt]
        & \le \mathbb{E}[\|x_k - \bar{x}_k\|^2] - 2 \alpha_k \mathbb{E}[ (f(x_{k+1}) - f(\bar{x}_k)) ]\nonumber\\[-1pt]
        & \quad - \frac{\rho \beta (1 -\beta)}{\blue{4} m c^2 B_h^2}\mathbb{E}[ \dist^2 (x_{k+1}, \mathcal{X})] + \alpha_k^2 B^2.
    \end{align*}
    
    \noindent Summing this from $0$ to $k-1$, and noticing that since $\alpha_k$ is a nonincreasing sequence, thus using $\alpha_k/ \alpha_0 \le 1$, we get:\vspace{-0.2cm}
    \begin{align*} 
        & \mathbb{E}[\|x_{k} - \bar{x}_{k}\|^2] \\[-3pt]
        & \le \|x_0 - \bar{x}_0\|^2 - 2 \sum_{t=0}^{k-1} \alpha_t \mathbb{E}[(f(x_{t+1}) - f(\bar{x}_{t}))]  \\[-2pt]
        & \quad  \!-\!  \frac{\rho \beta (1 - \beta)}{\blue{4} \alpha_0 m c^2 B_h^2} \sum_{t=0}^{k-1}  \alpha_t \mathbb{E}[\dist^2 (x_{t+1}, \mathcal{X})] \!+\! B^2 \sum_{t=0}^{k-1} \alpha_t^2.
    \end{align*}
    
    \noindent Now, divide the whole inequality by $S_k$ and from the definition of average sequences, linearity of the expectation operator and using convexity of the norm, and of the function $f$, we get:  
    \vspace{-0.2cm}
    \begin{align*} 
        & 2 \mathbb{E}[(f(\hat{x}_{k}) - f(\hat{x}^*_{k}))] + \frac{\rho \beta (1 - \beta)}{\blue{4}\alpha_0 m c^2 B_h^2} \mathbb{E}[\dist^2 (\hat{x}_{k}, \mathcal{X})]\\[-2pt]
        & \le \frac{\|x_0 - \bar{x}_0\|^2}{S_k} + B^2 \frac{\sum_{t=1}^{k} \alpha_t^2}{S_k}.
    \end{align*}
    Thus, we have the following rates in expectation for the average sequence in terms of optimality:\vspace{-0.15cm}
    \vspace{-0.2cm}
    \[ \mathbb{E}[(f(\hat{x}_{k}) - f(\hat{x}^*_{k}))] \le \frac{\|x_0 - \bar{x}_0\|^2}{2 S_k} + B^2 \frac{\sum_{t=1}^{k} \alpha_t^2}{2 S_k} . \] 
and feasibility violation:
\vspace{-0.2cm}
\begin{align*}
    & \mathbb{E}[\dist^2 (\hat{x}_{k}, \mathcal{X})] \\[-4pt]
    & \le \frac{\blue{4} \alpha_0 m c^2 B_h^2}{ \rho \beta (1 - \beta)} \left( \frac{\|x_0 - \bar{x}_0\|^2}{S_k} + B^2 \frac{\sum_{t=1}^{k} \alpha_t^2}{  S_k} \right).
\end{align*} 
Thus, we obtain the claimed results.
\end{proof}

\medskip

\noindent Now, Theorem \ref{Th:convex_case} yields (sublinear) convergence rates for SHAM iterates under convex
objective if the non-increasing stepsize $\alpha_k$ satisfies e.g., the conditions:  $\sum_{t=0}^{\infty}\alpha_t = \infty$ and $\sum_{t=0}^{\infty}\alpha_t^2 < \infty$ or $\sum_{t=0}^{k-1}\alpha_t^2 < \mathcal{O}(\ln(k+1))$ for all $k \geq 1$. Let us now  present two different choices for the stepsize $\alpha_k$:

\noindent \textbf{Choice} $(1)$:  Consider  $\alpha_k = \frac{\alpha_0}{\sqrt{k+2} \ln (k+2)}, \forall k \ge 1$, with $\alpha_0 \in \left(0, \frac{1}{L_f}\right]$.  Note that this choice will yield:
\vspace{-0.2cm}
    \begin{align*}
        & \sum_{t=1}^{k+1} \alpha_t \ge \frac{\alpha_0(k+1)}{\sqrt{k+3} \ln (k+3)}  \;\; \text{and } \;\; \sum_{t=1}^{k+1} \alpha_t^2 \le \frac{\alpha_0^2}{\ln(3)}.
    \end{align*}
    Thus, from Theorem \ref{Th:convex_case}, we obtain:
    \vspace{-0.2cm}
    \begin{align*}
        & \mathbb{E}[(f(\hat{x}_{k}) - f(\hat{x}^*_{k}))] \le \mathcal{O}\left(\frac{\ln (k+3)}{\sqrt{k+1}} \right),  \\[-3pt]
        & \mathbb{E}[\dist^2 (\hat{x}_{k}, \mathcal{X})] \le \mathcal{O}\left(\frac{\ln (k+3)}{\sqrt{k+1}} \right).
    \end{align*}


\noindent \textbf{Choice} $(2)$: Consider  $\alpha_k = \frac{\alpha_0}{\sqrt{k}}, \forall k \ge 1$, with $\alpha_0 \in \left(0, \frac{1}{L_f}\right]$.  This choice will give us:
\vspace{-0.19cm}
    \begin{align*}
        & \sum_{t=1}^{k+1} \alpha_t \ge \alpha_0 \sqrt{k+1},  \;\; \text{and } \sum_{t=1}^{k+1} \alpha_t^2 \le \mathcal{O} (\alpha_0^2 \ln (k+1)).
    \end{align*}
    Hence, from Theorem \ref{Th:convex_case}, we get:
    \vspace{-0.2cm}
    \begin{align*}
        & \mathbb{E}[(f(\hat{x}_{k}) - f(\hat{x}^*_{k}))] \le \mathcal{O}\left(\frac{1}{\sqrt{k+1}} + \frac{\ln (k+1)}{\sqrt{k+1}} \right), \\[-3pt]
        & \mathbb{E}[\dist^2 (\hat{x}_{k}, \mathcal{X})] \le \mathcal{O}\left(\frac{1}{\sqrt{k+1}} + \frac{\ln (k+1)}{\sqrt{k+1}} \right).
    \end{align*} 

\noindent Note that the second choice of $\alpha_k$, i.e., $\alpha_k = \frac{\alpha_0}{\sqrt{k}}$ provides better convergence rates (if we ignore the logarithmic term) for SHAM algorithm, which is $\mathcal{O} \left(  \frac{1}{\sqrt{k}}\right)$. To the best of our knowledge, these are the best possible rates for stochastic subgradient methods (see \cite{NecSin:22}). 

\vspace{-0.3cm}
\subsection{Convergence rates  for strongly convex objective}
\noindent In this section, we assume that all the assumptions made in the previous sections are true for problem \eqref{eq:prob} and give convergence rates for SHAM. First, we give a recurrence depending on $k_0 = \lfloor \frac{2L_f}{\mu} - 1\rfloor$, which gives rise to a switching stepsize strategy for the algorithm SHAM from constant stepsize to nonincreasing stepsize (here $x^*$ is the unique optimum of \eqref{eq:prob}). 
\vspace{-0.2cm}
\begin{lemma}\label{lemma_main_rec_strconv}
	\blue{Under the conditions of Theorem \ref{lemma_common}} with strongly convex $f$, i.e., $\mu>0$, define $k_0 = \lfloor \frac{2L_f}{\mu} - 1 \rfloor$,  $\theta_{\mu,L_f} \!=\! \left(1 \!-\! \frac{\mu}{L_f} \right)$ and $\alpha_k = \min \left( \frac{1}{L_f}, \frac{2}{\mu (k+1)}\right)$, we have the following recurrence true: 
 \vspace{-0.1cm}
	\begin{align}
        & \mathbb{E}[\|x_{k_0 + 1} - x^*\|^2] \\
        & \le \begin{cases}
            \frac{ B^2}{L_f^2}, \hspace{5.2 cm} \text{ if } \theta_{\mu, L_f} \le 0,\\
            \theta_{\mu,L_f}^{k_0+1} \mathbb{E}[\|x_0 - x^*\|^2 ] + (1 - \theta_{\mu,L_f}^{k_0+1}) \frac{ B^2}{\mu L_f}, \text{ if, } \theta_{\mu,L_f} > 0,
        \end{cases} \nonumber\\
		&(k+1)^2\mathbb{E}[\|x_{k+1} - x^*\|^2] \le k^2\mathbb{E}[\|x_k - x^*\|^2]   \label{eq:main_rec_strconv2}\\
        & \quad - \frac{4 (k+1)}{\mu}\mathbb{E}[(f(x_{k+1}) - f(x^*))] + \frac{4B^2}{\mu^2}\nonumber\\
        & \;\;\; -  \frac{\rho \beta (1 - \beta) (k\!+\!1)^2}{\blue{4} m c^2 B_h^2}\mathbb{E}[\dist^2 (x_{k+1}, \mathcal{X})], \;\forall k > k_0.\nonumber
	\end{align}
\end{lemma}

\begin{proof} 
See Appendix. 
\end{proof}

\noindent Now, for $k \geq k_0+1$, let us define the sum:
\vspace{-0.3cm}
\begin{align*}
	\bar{S}_k = \sum_{t=k_0+1}^{k} (t+1)^2 \sim \mathcal{O} (k^3 + k_0^2k + k^2k_0).
\end{align*}
and the corresponding average sequences:
\vspace{-0.3cm}
\begin{align*}
	&\hat{x}_k = \frac{1}{\bar{S}_k} \sum_{t=k_0+1}^{k} (t+1)^2 x_{t+1}, \\[-2pt]
	&  \hat{w}_k = \frac{\sum_{t=k_0+1}^{k} (t+1)^2 \Pi_{\mathcal X}(x_{t+1})}{S_k} \in \mathcal{X}. 
\end{align*}
\vspace{0.03cm}

\noindent The following theorem provides the convergence rates of Algorithm \ref{alg:alg1} under Assumption \ref{assumption1}$(i)$.

\begin{theorem}\label{Th:strconvex_case}
	\blue{Under the conditions of Lemma \ref{lemma_main_rec_strconv}}, we have the following convergence rates for the average sequence $\hat{x}_{k}$ in terms of optimality and feasibility violation:
 \vspace{-0.05cm}
	\begin{align*}
		& \mathbb{E}[(f(\hat{x}_{k}) - f^*)] \le \mathcal{O} \left( \frac{B^2}{\mu^2 (k-k_0)} \right),\\[-2pt]
		& \mathbb{E}[\dist^2 (\hat{x}_{k}, \mathcal{X})] \le \mathcal{O} \left( \frac{\blue{4} m c^2 B_h^2B^2}{\rho \beta (1 - \beta)\mu^2 (k^2 \!+\! k k_0 \!+\! k_0^2)} \right).
	\end{align*}
\end{theorem}
\medskip
\begin{proof}
	From Lemma \ref{lemma_main_rec_strconv}, for any $k > k_0$, we have \eqref{eq:main_rec_strconv2}. Summing it from $k_0 + 1$ to $k$, we get:
 \vspace{-0.2cm}
\begin{align*} 
	& (k+1)^2 \mathbb{E}[\|x_{k+1} - x^*\|^2] \le (k_0+1)^2 \|x_{k_0 + 1} - x^*\|^2 \\[-2pt]
	& \quad - \frac{4}{\mu} \sum_{t=k_0 + 1}^{k} (t+1) \mathbb{E}[(f(x_{t+1}) - f(x^*))] \\[-3pt]
	& \quad  -  \frac{\rho \beta (1 - \beta)}{\blue{4}  m c^2 B_h^2} \sum_{t=k_0 + 1}^{k} (t+1)^2 \mathbb{E}[\dist^2 (x_{t+1}, \mathcal{X})] \\[-4pt]
    & \quad + \frac{4B^2}{\mu^2} (k - k_0).
\end{align*}


\noindent By linearity of the expectation operator and using convexity of the norm and of the function $f$, we obtain:
\vspace{-0.1cm}
\begin{align*} 
	& \frac{4\bar{S}_k}{\mu (k\!+\! 1)} \mathbb{E}[(f(\hat{x}_{k}) \!-\! f^*)] \!+\! \frac{\rho \beta (1 - \beta) \bar{S}_k}{\blue{4} m c^2 B_h^2}\mathbb{E}[\dist^2 (\hat{x}_{k}, \mathcal{X})]\\[-3pt]
	& \le (k_0+1)^2 \|x_{k_0 + 1} - x^*\|^2  + \frac{4B^2}{\mu^2} (k - k_0).
\end{align*}

\noindent Thus, we have the following rates in expectation for the average sequence in terms of optimality:
\vspace{-0.1cm}
\begin{align*}
    & \mathbb{E}[(f(\hat{x}_{k}) - f^*)] \\[-3pt]
    & \le \frac{\mu(k+1)}{4\bar{S}_k}\left( (k_0+1)^2 \|x_{k_0 + 1} - x^*\|^2 + \frac{4B^2}{\mu^2} (k - k_0)\right) \\[-3pt]
    & \approx \mathcal{O} \left( \frac{B^2}{\mu^2 (k-k_0)} \right),
\end{align*} 
and feasibility violation:
\vspace{-0.1cm}
\begin{align*}
    & \mathbb{E}[\dist^2 (\hat{x}_{k}, \mathcal{X})] 
    \approx \mathcal{O} \left( \frac{\blue{4} m c^2 B_h^2B^2}{\rho \beta (1 - \beta) \mu^2 (k^2 \!+\! k k_0 \!+\! k_0^2)} \right).
\end{align*} 
This concludes the proof of the theorem.
\end{proof}

\noindent Since we are using the $\mu$-strongly convex condition here, the  convergence rates are now $\mathcal{O} \left(\frac{1}{k}\right)$ for SHAM algorithm.

\vspace{-0.3cm}
\section{Numerical results}

\noindent In this section, we consider a problem characterized by a quadratic objective function and second order cone (SOC) constraints, presented as follows: 
\begin{equation}\label{eq:numerical_experiment}
	\begin{array}{rl}
	& \min\limits_{x \red{\in \mathcal{Y}}\in \mathbb{R}^n} \; \frac{1}{2}x^TQ_f x + q_f^T x\\
	& \text{subject to } \;  \|Q_i x + a_i\| \le q_i^T x + b_i \;\; \forall i = 1:m,
	\end{array}
\end{equation}

\noindent where $Q_f \in \mathbb{R}^{n\times n}, q_f \in \mathbb{R}^{n}, Q_i \in \mathbb{R}^{n_i\times n}, a_i \in \mathbb{R}^{n_i}, q_i \in \mathbb{R}^{n}$, $b_i\in \mathbb{R}_+$ \red{and compact set  $\mathcal{Y} = [-10^3, \; 10^3]$} (hence, the feasible set contains an open set around origin) represent the problem parameters, with $x$ as the optimization variable. Hence, this problem aligns with all the assumptions of problem \eqref{eq:prob}. Depending on whether $Q_f$ is positive semidefinite or positive definite, \eqref{eq:numerical_experiment} possesses either a convex or $\mu$-strongly convex objective function, respectively. Furthermore, problem \eqref{eq:numerical_experiment} exhibits remarkable generality, as it encompasses convex quadratic programs (QPs), quadratically constrained quadratic programs (QCQPs), and numerous other nonlinear convex optimization problems, as elucidated in \cite{LobLeb:98}. Moreover, SOC constraints are general over quadratic constraints in convex settings (see Section $A.2.3$ in \cite{Wei:20}).  From a practical perspective, problem \eqref{eq:numerical_experiment} finds relevance in diverse fields. For instance, \cite{BoyHan:97} (equation $11$) illustrates its applicability in robust optimal control problems. In the realm of signal processing, it emerges as a powerful tool for robust beamformer design, facilitating the minimization of transmitted power while adhering to received signal-to-noise ratio constraints, as showcased in \cite{SidLuo:06}.

\noindent Practical implementations are conducted using MATLAB R2023b on a laptop equipped with an i5 CPU operating at 2.1 GHz and 16 GB of RAM. The problem parameters are generated as random data within MATLAB, with the added constraint that $n_i < n$  for all $i \in [m]$. 

\vspace{0.1cm}

\noindent Table \ref{table} below provides a comprehensive overview of CPU times in seconds, encompassing minimum, average, and maximum values derived from $5$ runs. These times are measured for both,  CVX solver \cite{GraBoy:13} and  SHAM algorithm. \blue{The best average time taken by SHAM to solve \eqref{eq:numerical_experiment} is written in bold}. Two distinct choices for $\tilde{x}_k$ are explored, namely, $x_k \, (\gamma = 0)$ and  $v_k \, (\gamma = 1)$, with the stepsize $\alpha_k$ defined as $\alpha_k = \frac{\alpha_0}{(\ln(k+1)\sqrt{k+1})}$, where $\alpha_0 = \frac{1}{L_f}$ (convex case with $\mu = 0$), and $\alpha_k = \min\left(\frac{1}{L_f}, \frac{2}{(\mu (k+1))}\right)$ (strongly convex case with $\mu > 0$). Furthermore, we set $\beta = 0.96$ and each index $j_k$ is chosen uniformly at random, considering an epoch of the algorithm as complete when the number of iterations matches the number of constraints. Our stopping criteria encompass $\|\max(0, h(x))\|^2 \le 10^{-2}$ and $|f(x) - f^*| \leq 10^{-2}$ (with $f^*$ computed via CVX solution) or $\max (\|x_{k+1} - x_k\|^2, \ldots, \|x_{k-M+1} - x_{k-M}\|^2) \le 10^{-3}$, where we take $M=10$ (when CVX does not solve in $6$ hours).


\begin{figure}[ht]
	\centering 
	\includegraphics[height=4.3cm, width=4.35cm]{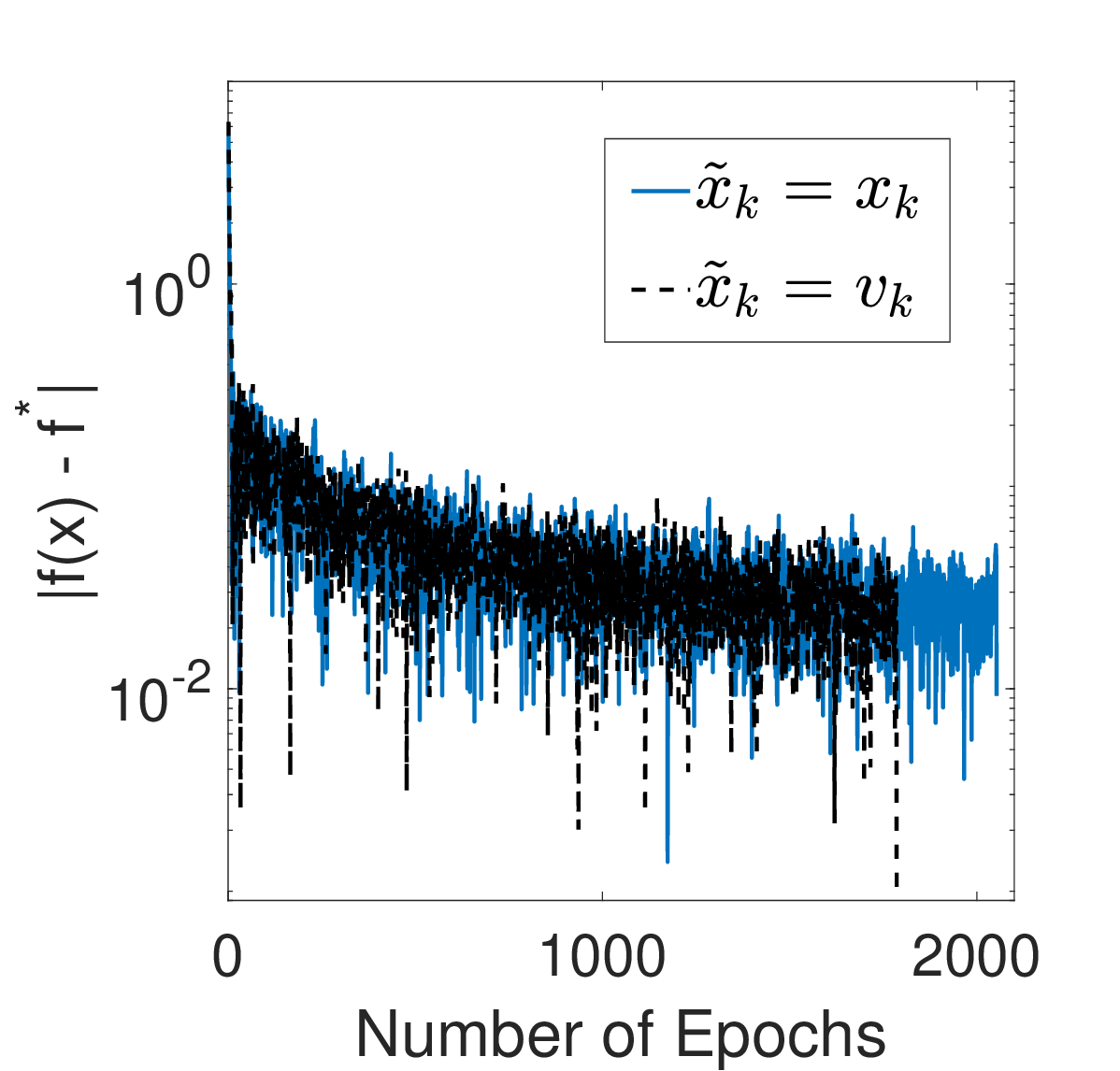}
	\includegraphics[height=4.3cm, width=4.35cm]{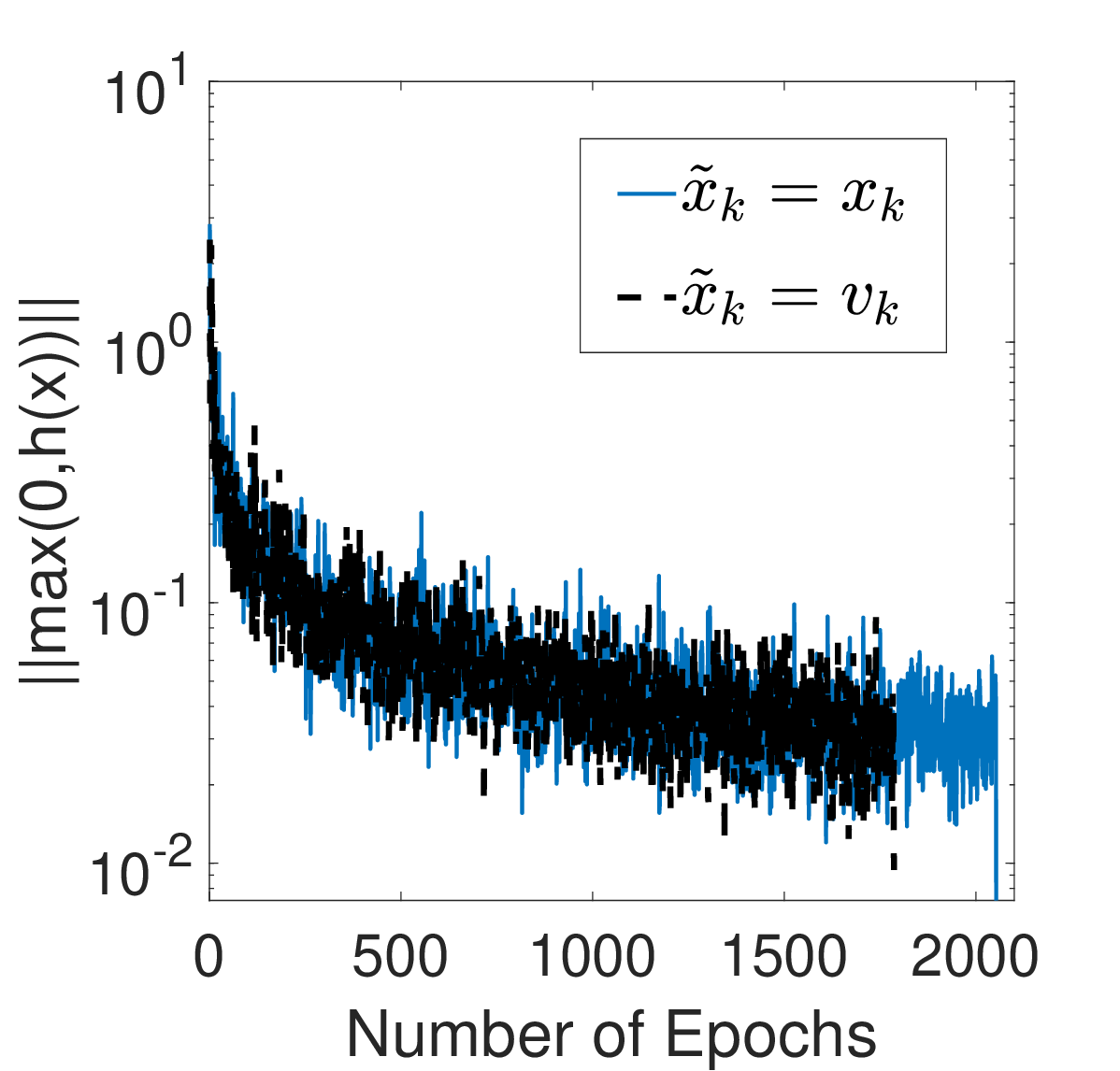}
        \includegraphics[height=4.3cm, width=4.35cm]{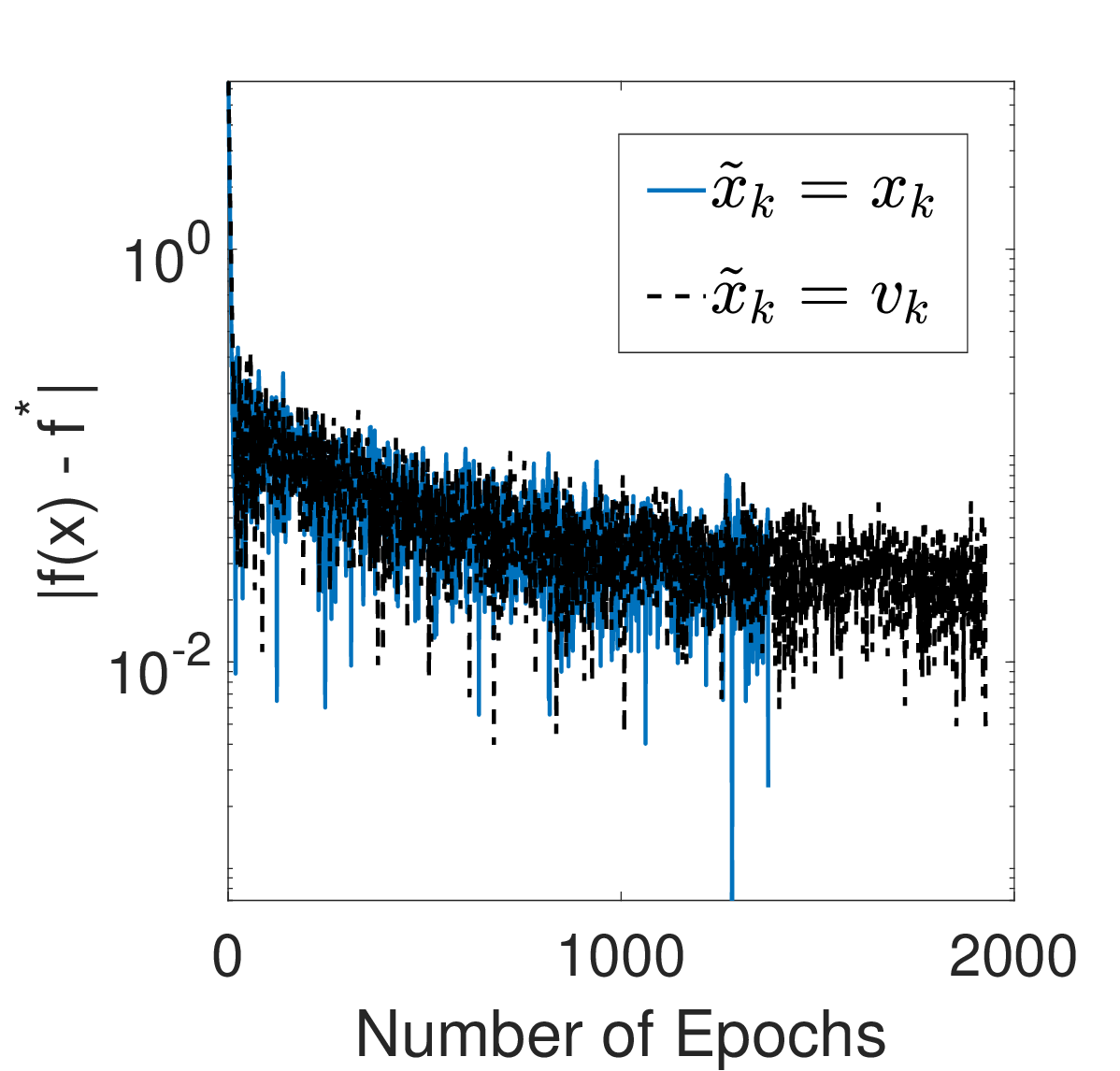}
	\includegraphics[height=4.3cm, width=4.35cm]{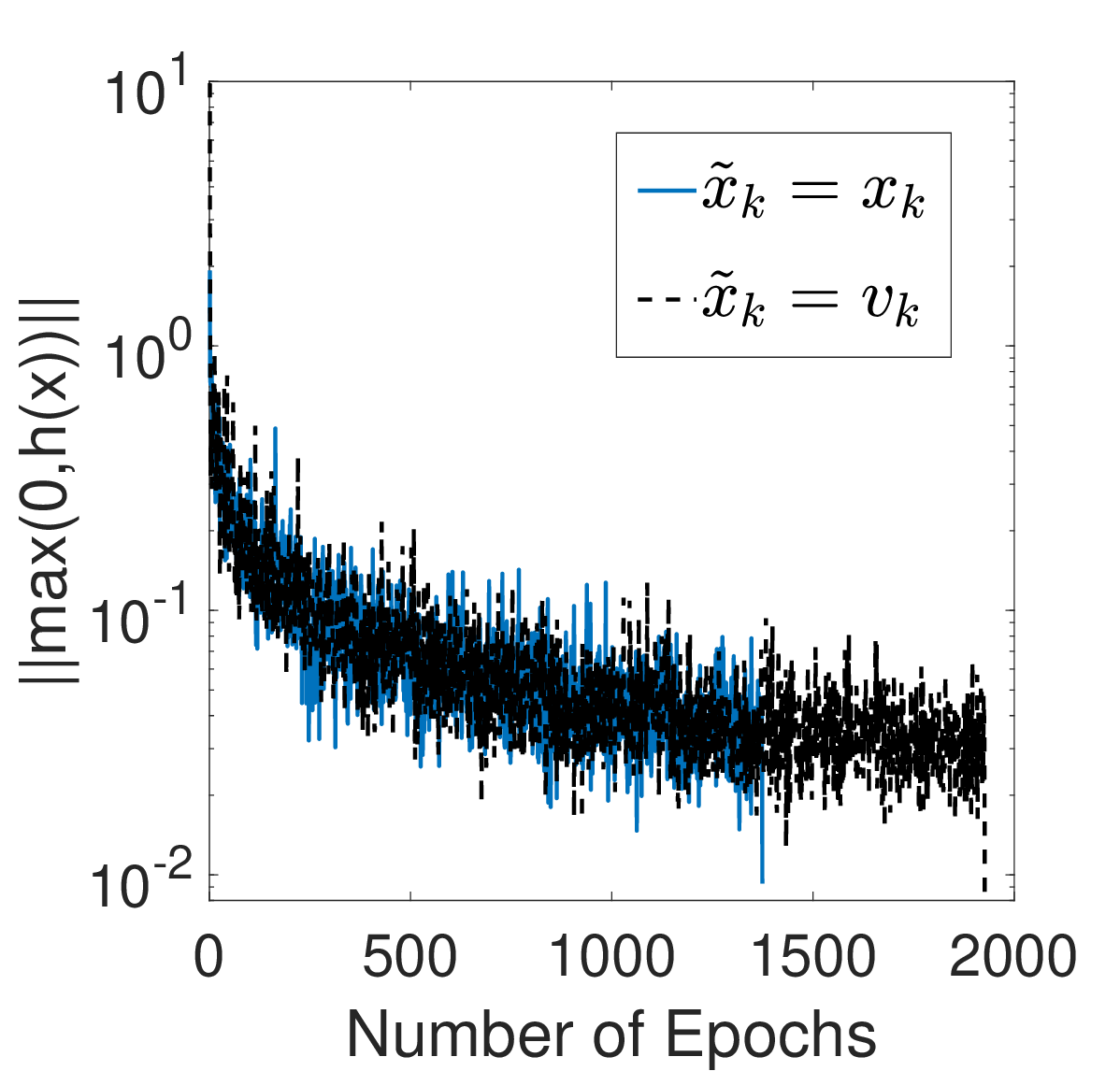}
	\caption{Behaviour of SHAM in terms of optimality (left) and feasibility (right) running two times (top and bottom) on the same problem with $ n=100,  m = 100,  \mu = 0$, and two choices for  $\tilde{x}_k$.}
	\label{fig_alpha_conv}
 \vspace{-0.5cm}
\end{figure}

\begin{table}[ht]\centering 
\caption{Performance of SHAM versus  CVX.}
\label{table}
\begin{tabular}{|c|clll|ll|}
\hline
\multirow{3}{*}{Data-size}                               & \multicolumn{4}{c|}{SHAM}                                                                                                                                                      & \multicolumn{2}{c|}{\multirow{2}{*}{CVX}}                 \\ \cline{2-5}
                                                         & \multicolumn{2}{c|}{$\tilde{x}_k = x_k$}                       & \multicolumn{2}{c|}{$\tilde{x}_k = v_k$}                                                                      & \multicolumn{2}{c|}{}                                     \\ \cline{2-7} 
                                                         & \multicolumn{1}{c|}{$\mu = 0$} & \multicolumn{1}{c|}{$\mu> 0$} & \multicolumn{1}{c|}{$\mu = 0$} & \multicolumn{1}{c|}{$\mu> 0$}                                                & \multicolumn{1}{l|}{$\mu = 0$} & $\mu>0$                  \\ \hline
\begin{tabular}[c]{@{}c@{}}n = $10^2$\\ m = $10^2$\end{tabular} & \multicolumn{1}{c|}{\begin{tabular}[c]{@{}c@{}}8.4\\ \textbf{11.3}\\ 13.8\end{tabular}}       & \multicolumn{1}{c|}{\begin{tabular}[c]{@{}c@{}}4.5\\ 11.6\\ 18.4\end{tabular}}      & \multicolumn{1}{c|}{\begin{tabular}[c]{@{}c@{}}11.2\\ 12.3\\ 13.2\end{tabular}}       & \multicolumn{1}{c|}{\begin{tabular}[c]{@{}c@{}}4.7\\ \textbf{10.8}\\ 17.7\end{tabular}} & \multicolumn{1}{c|}{2.4}       & \multicolumn{1}{c|}{1.9} \\ \hline
\begin{tabular}[c]{@{}c@{}}n = $10^2$\\ m = $10^3$\end{tabular} & \multicolumn{1}{c|}{\begin{tabular}[c]{@{}c@{}} 75.3\\ 86.4\\102.4 \end{tabular}}       & \multicolumn{1}{c|}{\begin{tabular}[c]{@{}c@{}} 75.3\\ \textbf{105.9}\\156.4 \end{tabular}}      & \multicolumn{1}{c|}{\begin{tabular}[c]{@{}c@{}} 52.2\\ \textbf{86.1}\\ 108.7\end{tabular}}       & \multicolumn{1}{c|}{\begin{tabular}[c]{@{}c@{}} 124.5\\ 142.4 \\ 185.2\end{tabular}} & \multicolumn{1}{c|}{95.5}       & \multicolumn{1}{c|}{101.6} \\ \hline
\begin{tabular}[c]{@{}c@{}}n = $10^2$\\ m = $10^4$\end{tabular} & \multicolumn{1}{c|}{\begin{tabular}[c]{@{}c@{}}1898.2\\ \textbf{2326.1}\\2694.6 \end{tabular}}       & \multicolumn{1}{c|}{\begin{tabular}[c]{@{}c@{}}181.9\\ \textbf{541.3}\\ 1110.4\end{tabular}}      & \multicolumn{1}{c|}{\begin{tabular}[c]{@{}c@{}}2217.9\\ 2503.6\\ 2755.5\end{tabular}}       & \multicolumn{1}{c|}{\begin{tabular}[c]{@{}c@{}}264.4\\ 716.6 \\ 1095.4\end{tabular}} & \multicolumn{1}{c|}{7433.8}       & \multicolumn{1}{c|}{7344.9} \\ \hline
\begin{tabular}[c]{@{}c@{}}n = $10^3$\\ m = $10^2$\end{tabular} & \multicolumn{1}{c|}{\begin{tabular}[c]{@{}c@{}}181.6\\ 226.9\\ 275.5\end{tabular}}       & \multicolumn{1}{c|}{\begin{tabular}[c]{@{}c@{}}88.5\\ \textbf{117.7}\\ 155.1\end{tabular}}      & \multicolumn{1}{c|}{\begin{tabular}[c]{@{}c@{}}177.6\\ \textbf{212.3}\\241.5 \end{tabular}}       & \multicolumn{1}{c|}{\begin{tabular}[c]{@{}c@{}}118.1\\ 153.6\\ 191.5\end{tabular}} & \multicolumn{1}{c|}{2260.3}       & \multicolumn{1}{c|}{1837.2} \\ \hline
\begin{tabular}[c]{@{}c@{}}n = $10^3$\\ m = $10^3$\end{tabular} & \multicolumn{1}{c|}{\begin{tabular}[c]{@{}c@{}}2519.2\\ \textbf{2926.9}\\3271.1 \end{tabular}}       & \multicolumn{1}{c|}{\begin{tabular}[c]{@{}c@{}}395.3\\ \textbf{426.3}\\458.9 \end{tabular}}      & \multicolumn{1}{c|}{\begin{tabular}[c]{@{}c@{}}2547.5\\ 3567.1\\ 4436.3\end{tabular}}       & \multicolumn{1}{c|}{\begin{tabular}[c]{@{}c@{}}2124.1\\ 2560.5\\2884.1 \end{tabular}} & \multicolumn{1}{c|}{*}       & \multicolumn{1}{c|}{*} \\ \hline
\end{tabular}
\vspace{-0.7cm}
\end{table}

\noindent From the two choices of $\tilde{x}_k$, $\tilde{x}_k = x_k$ and $\tilde{x}_k = v_k$ in Table I,  intuitively, one can say that the choice $\tilde{x}_k = v_k$ should be better, but in practice $\tilde{x}_k = x_k$ can be much faster. However, the choice $\tilde{x}_k = v_k$ is more robust in the sense that it has less variation in maximum and minimum CPU timings  (see also  Figure \ref{fig_alpha_conv}). Comparing the performance of SHAM against CVX, our findings showcase that SHAM consistently outperforms CVX by a factor ranging from $3$ to $10$, especially when dealing with high dimensions. Notably, our analysis also highlights a significant advantage of SHAM in distinguishing between scenarios  $\mu = 0$ and $\mu > 0$, a distinction not made by CVX solver, which relies on interior point methods.

\vspace{-0.3cm}





\begin{appendix}
    \noindent(1) \textit{Proof of Lemma \ref{lemma_lin_reg}:} If $x \in \mathcal{X}\blue{\subseteq \mathcal{Y}}$, the result holds trivially. So let us assume that $x \notin \mathcal{X}$. Define the index set 
    $\mathcal{I}(x) = \{ j \in [m] : h_j(x) > 0 \}$.
    Since $x \in \blue{\mathcal{Y}} \backslash \mathcal{X}$,  $\mathcal{I}(x)$ is nonempty. From \eqref{eq:lin_reg}:\vspace{-0.2cm}
    \begin{align}\label{eq:1}
        \dist (x, \mathcal{X}) \le c \max_{j\in [m]} \left[ (h_j(x))_+ \right] = c h_{j^*}(x),
    \end{align}
    where $j^* \in \argmax_j \{ h_j(x) : j \in [m] \}$ and $c$ is independent of $x$. Clearly $j^* \in \mathcal{I}(x)$ and  $h_{j^*}(x) = \max_{j \in \mathcal{I}(x)} \{h_j(x)\}$. Fix  $\Delta_{j^*}^x \in \partial h_{j^*} (x)$. If $\Delta_{j^*}^x \neq 0$, from \eqref{eq:half_space} and the fact that $\Pi_{\mathcal{L}(h_{j^*}; x; \Delta_{j^*}^x)} (x) \in \mathcal{L}(h_{j^*};x;\Delta_{j^*}^x)$, we get:
    \vspace{-0.2cm}
    \begin{align*}
        h_{j^*} \!(x) \!& \le\! \!-\langle \Delta_{j^*}^x\!, \! \Pi_{\mathcal{L}(h_{j^*}\!;x;\Delta_{j^*}^x\!)} (x) \!-\! x \rangle \\[-2pt]
        & \le\! B_h \dist(x, \!\mathcal{L}(h_{j^*}\!;x;\Delta_{j^*}^x\!)\!),
    \end{align*} 
    where the last inequality follows from the Cauchy-Schwartz inequality and Assumption \ref{assumption2}. 
     Minimizing the right-hand side over $\Delta_{j^*}^x$
     and using  \eqref{eq:1}, we get:
     \vspace{-0.25cm}
     \begin{align*}
         & \dist (x, \mathcal{X}) \le c B_h \min_{\Delta_{j^*}^x \in \partial h_{j^*}(x)} \dist(x, \mathcal{L}(h_{j^*};x;\Delta_{j^*}^x))\\[-4pt]
         & \le c B_h \max_{j\in[m]} \min_{\Delta_{j}^x \in \partial h_{j}(x)} \dist(x, \mathcal{L}(h_{j};x;\Delta_{j}^x)).
     \end{align*}
     If $\Delta_{j^*}^x = {0}$, from convexity, $x$ is the minimizer of $h_{j^*}$, and $h_{j^*} (x) = \min_{y \in \mathbb{R}^n} h_{j^*} (y) \le 0,$
     where the inequality follows from the assumption that $\mathcal{X}$ is non-empty. However, by the definition of $j^*$, $h_{j^*} (x)>0$. Hence it is a contradiction, thus it is impossible to have $\Delta_{j^*}^x = {0}$.\\
     
     (2) \textit{Proof of Lemma \ref{lemma_smooth_f}:} From the smoothness property of the function $f$, we have:
     \vspace{-0.3cm}
    \begin{align*}
	& 2\langle \blue{u_k} - x_{k+1}, x_{k+1} - x_k \rangle  \\[-3pt]
        & \overset{\eqref{eq:alg2step1}}{=} 2 \langle x_k - x_{k+1} - \alpha_k \nabla f(x_k), x_{k+1} - x_k \rangle\\[-3pt]
	& \overset{\eqref{eq:smooth_f}}{\le} 
		2 \alpha_k ( f(x_k) - f(x_{k+1})) + (\alpha_k L_f - 2)\|x_{k+1} - x_k\|^2.
    \end{align*}
    Hence we get the desired result.\\
    \vspace{-0.3cm}
    
    (3) \textit{Proof of Lemma \ref{lemma_non_smooth_h2}:} Following \cite{YanYue:22} (Lemma 4), we have:\vspace{-0.1cm}
    \begin{align*}
        & \mathbb{E}[ \|\tilde{x}_k - \Pi_{\mathcal{L}_{j_k}} (\tilde{x}_k)\|^2 |\mathcal{F}_{[k]}]\\[-3pt]
        & \ge \frac{\rho}{m} \min_{\Delta_{j}^{\tilde{x}_k} \in \partial h_j(\tilde{x}_k)} \dist^2 (\tilde{x}_k,\mathcal{L} (h_j; \tilde{x}_k; \Delta_{j}^{\tilde{x}_k})).
    \end{align*}
    Noticing that $\tilde{x}_k \in \blue{\mathcal{Y}}$, from Lemma \ref{lemma_lin_reg}, we get:\vspace{-0.28cm}
    \begin{align*}
        & \mathbb{E}[ \|\tilde{x}_k - \Pi_{\mathcal{L}_{j_k}} (\tilde{x}_k)\|^2 |\mathcal{F}_{[k]}]
        \overset{\eqref{eq:lin_reg2}}{\ge}\frac{\rho}{m c^2 B_h^2} \dist^2 (\tilde{x}_k, \mathcal{X}).
    \end{align*}
    After rearranging the terms, we obtain:\vspace{-0.2cm}
    \begin{align*}
        & \dist^2 (\tilde{x}_k, \mathcal{X}) \le \frac{m c^2 B_h^2}{\rho} \mathbb{E}\left[ \|\tilde{x}_k - \Pi_{\mathcal{L}_{j_k}} (\tilde{x}_k)\|^2 |\mathcal{F}_{[k]}\right]\\[-3pt]
        & \le \frac{2 m c^2 B_h^2}{\rho} \left(\mathbb{E}\left[\|v_k - \Pi_{\mathcal{L}_{j_k}} (v_k)\|^2 |\mathcal{F}_{[k]}\right] + \|v_k - \tilde{x}_k\|^2\right),
    \end{align*}
    where in the last inequality we use \eqref{eq:proj_property} and \eqref{eq:ineq2}.\\

    (4) \textit{Proof of Lemma \ref{lemma_dist}:} For any $\tilde{x}_k \in \blue{\mathcal{Y}}$ and given $\Delta_{j_k}^{\tilde{x}_k} \in \partial h_{j_k}(\tilde{x}_k)$, consider the following two cases.\\
    \noindent \textit{Case (i):} When $\Delta_{j_k}^{\tilde{x}_k} = 0$, from \eqref{eq:algstep2} we have $\red{z_k} = v_k$. Thus the result holds automatically after using \eqref{eq:proj_property} and \eqref{eq:proj_property2}.\\
    \noindent \textit{Case (ii):} When $\Delta_{j_k}^{\tilde{x}_k} \neq 0$, \red{using $\Pi_\mathcal{X} (v_k) \in \mathcal{X} \subseteq \mathcal{Y}$}, we have:	
    \vspace{-0.4cm}
    \begin{align*}
        & \dist^2 (x_{k+1}, \mathcal{X}) 
        \overset{\eqref{eq:proj_property}}{\le}  \|x_{k+1} - \Pi_\mathcal{X} (v_k) \|^2 \red{\overset{\eqref{eq:proj_property2}}{\le}  \|z_k - \Pi_\mathcal{X} (v_k) \|^2} \\[-4pt]
        & \overset{\eqref{eq:algstep}}{=} \|v_k - \Pi_\mathcal{X} (v_k)\|^2 + \beta^2 \frac{ (l_{h_{j_k}}(v_k; \tilde{x}_k))^2_+ }{ \|\Delta_{j_k}^{\tilde{x}_k}\|^2} \nonumber\\[-4pt]
        & \quad - 2\beta \frac{ (l_{h_{j_k}}(v_k; \tilde{x}_k))_+}{\|\Delta_{j_k}^{\tilde{x}_k}\|^2} \langle \Delta_{j_k}^{\tilde{x}_k}, v_k - \Pi_\mathcal{X} (v_k) \rangle \\[-4pt]
        & \le \|v_k - \Pi_\mathcal{X} (v_k)\|^2 + \beta^2 \frac{ (l_{h_{j_k}}(v_k; \tilde{x}_k))^2_+ }{ \|\Delta_{j_k}^{\tilde{x}_k}\|^2}\nonumber \\[-4pt]
        & \quad  - 2\beta \frac{ (l_{h_{j_k}}(v_k; \tilde{x}_k))_+}{\|\Delta_{j_k}^{\tilde{x}_k}\|^2} (h_{j_k}(\tilde{x}_k) - h_{j_k}(\Pi_\mathcal{X} (v_k)))\\[-4pt]
        & \quad - 2\beta \frac{ (l_{h_{j_k}}(v_k; \tilde{x}_k))_+}{\|\Delta_{j_k}^{\tilde{x}_k}\|^2} \langle \Delta_{j_k}^{\tilde{x}_k}, v_k - \tilde{x}_k \rangle\\[-4pt]
        & \le  \|v_k - \Pi_\mathcal{X} (v_k)\|^2 + \beta^2 \frac{ (l_{h_{j_k}}(v_k; \tilde{x}_k))^2_+ }{ \|\Delta_{j_k}^{\tilde{x}_k}\|^2}\\[-4pt]
        & \quad - 2\beta \frac{ (l_{h_{j_k}}(v_k; \tilde{x}_k))_+}{\|\Delta_{j_k}^{\tilde{x}_k}\|^2} l_{h_{j_k}}(v_k; \tilde{x}_k) \\[-4pt]
        & \overset{\eqref{eq:ineq3}}{=} \|v_k - \Pi_\mathcal{X} (v_k)\|^2 - \frac{\beta (2 - \beta)}{\|\Delta_{j_k}^{\tilde{x}_k}\|^2}(l_{h_{j_k}}(v_k; \tilde{x}_k))^2_+,
    \end{align*}
    where the third inequality follows from  convexity of  function $h_{j_k}$ and the fourth inequality uses that $h(\Pi_\mathcal{X} (v_k)) \le 0$. Further, noticing that $\beta \in (0,2)$, we obtain the statement.\\

\setlength{\belowdisplayskip}{0.03cm} 
\setlength{\belowdisplayshortskip}{0.03cm}

    (5) \textit{Proof of Lemma \ref{lemma_non_smooth_h}:} From  Lemma \ref{lemma_dist}, we have:
    \vspace{-0.1cm}
    \begin{align*}
        & \dist^2 (x_{k+1}, \mathcal{X}) \le \|v_k - \Pi_\mathcal{X} (v_k)\|^2 \overset{\eqref{eq:proj_property}}{\le} \|v_k - \Pi_\mathcal{X} (\tilde{x}_k) \|^2\\
        & \overset{\eqref{eq:ineq2}}{\le} 2 \|v_k - \tilde{x}_k\|^2 + 2 \|\tilde{x}_k - \Pi_\mathcal{X} (\tilde{x}_k) \|^2.
    \end{align*}
    Now, after taking expectation conditioned on $j_k$, and using the result from Lemma \ref{lemma_non_smooth_h2}, we get:
    \vspace{-0.1cm}
    \begin{align*}
        & \mathbb{E}[ \dist^2 (x_{k+1}, \mathcal{X})|\mathcal{F}_{[k]}] 
        \le \left( 2 + \frac{4 m c^2 B_h^2}{\rho}\right) \|v_k - \tilde{x}_k\|^2  \\
        & \quad + \frac{4 m c^2 B_h^2}{\rho}  \mathbb{E}\left[\|v_k - \Pi_{\mathcal{L}_k} (v_k)\|^2 |\mathcal{F}_{[k]}\right].
    \end{align*} 
    After rearranging the terms we get the required result. \\

    (6) \textit{Proof of Lemma \ref{lemma_main_rec_strconv}:} From Theorem \ref{lemma_common}, we have:\vspace{-0.1cm}
\begin{align} \label{eq:mid_way1}
	& \mathbb{E}[\|x_{k+1} - x^*\|^2] \\
	& \le (1 -\mu \alpha_k)\mathbb{E}[\|x_k - x^*\|^2] - 2 \alpha_k \mathbb{E}[(f(x_{k+1}) - f(x^*))]  \nonumber\\
    & \quad  -  \frac{\rho\beta (1 - \beta)}{\blue{4} m c^2 B_h^2}\mathbb{E}[\dist^2 (x_{k+1}, \mathcal{X})] + \alpha_k^2 B^2\nonumber.
\end{align}
For $k\le k_0$, we have $\alpha_k = \frac{1}{L_f}$, thus from \eqref{eq:mid_way1}, we get:\vspace{ - 0.1cm}
\begin{align*}
    & \mathbb{E}[\|x_{k+1} - x^*\|^2] \le \left(1 - \frac{\mu}{L_f} \right) \mathbb{E}[\|x_k \!-\! x^*\|^2 ] + \frac{ B^2 }{L_f^2}\\
    & \!\le\! \max \left(  \left(1 - \frac{\mu}{L_f} \right) \mathbb{E}[\|x_k \!-\! x^*\|^2 ] + \frac{ B^2 }{L_f^2}, \frac{ B^2 }{L_f^2}\right).
\end{align*}
Using the geometric sum formula and noticing that $\theta_{\mu,L_f} = \left(1 - \frac{\mu}{L_f} \right)$, we obtain the first statement. 
Next, for $k > k_0$, from the choice of $\alpha_k$, we know that $\alpha_k = \frac{2}{\mu(k+1)}$, thus from \eqref{eq:mid_way1}, we have:
\vspace{-0.1cm}
\begin{align*}
    & \mathbb{E}[\|x_{k+1} - x^*\|^2] \le \left(1 -\frac{2}{k+1}\right)\mathbb{E}[\|x_k - x^*\|^2]  \nonumber\\[-3pt]
    & \quad - \frac{4}{\mu(k+1)} \mathbb{E}[(f(x_{k+1}) - f(x^*))]  \\
    & \quad -  \frac{\rho \beta (1 - \beta)}{\blue{4} m c^2 B_h^2} \mathbb{E}[\dist^2 (x_{k+1}, \mathcal{X})] + \frac{4B^2}{\mu^2 (k+1)^2}\\
    & = \frac{k-1}{k+1}\mathbb{E}[\|x_k - x^*\|^2] - \frac{4}{\mu(k+1)} \mathbb{E}[(f(x_{k+1}) - f(x^*))]  \nonumber\\[-3pt]
    & \quad  -  \frac{\rho \beta (1 -\beta)}{\blue{4} m c^2 B_h^2} \mathbb{E}[\dist^2 (x_{k+1}, \mathcal{X})] + \frac{4B^2}{\mu^2 (k+1)^2} .
\end{align*}
\noindent Now, multiply the whole inequality by $(k+1)^2$, and using the fact that $k^2 - 1\le k^2$, we get:
\vspace{-0.1cm}
\begin{align*}
    & (k+1)^2\mathbb{E}[\|x_{k+1} - x^*\|^2] \le k^2 \mathbb{E}[\|x_k - x^*\|^2]\\
    & \quad - \frac{4(k+1)}{\mu} \mathbb{E}[(f(x_{k+1}) - f(x^*))]\\
    & \quad -  \frac{\rho \beta (1 - \beta) (k+1)^2}{\blue{4}m c^2 B_h^2} \mathbb{E}[\dist^2 (x_{k+1}, \mathcal{X})]  + \frac{4B^2}{\mu^2}.
\end{align*}

\noindent Hence we get the result \eqref{eq:main_rec_strconv2}.
\end{appendix}






\section*{References}

\begin{thebibliography}{00}\leftskip1pc
	

\vspace{-0.4cm}

    \bibitem{BerBac:22}
    E. Berthier, J. Carpentier, A. Rudi and F. Bach, "Infinite-Dimensional Sums-of-Squares for Optimal Control", {\it Conference on Decision and Control}, 577-582, 2022.

    \bibitem{BhaGra:04}
	C.  Bhattacharyya,  L.R. Grate, M.I. Jordan, L. El Ghaoui and S. Mian, "Robust sparse hyperplane classifiers: application to uncertain molecular profiling data,"  {\it Journal of Computational Biology}, vol. 11, no. 6, 1073-1089, 2004.


    \bibitem{BoyHan:97}
    S. Boyd, C. Crusius and A. Hansson, "Control applications of nonlinear convex programming," {\it Journal of Process Control}, vol. 8, no. 5-6, 313-324, 1998.

    \bibitem{BooLan:23}
    \red{D. Boob, Q. Deng and G. Lan, "Stochastic first-order methods for convex and nonconvex functional constrained optimization", \textit{Mathematical Programming}, vol. 197, 215–279, 2023.}





    \bibitem{GowRic:19}
	R. Gower, L. Nicolas, Q. Xun, S. Alibek, S. Egor and P. Richtarik, "SGD: general analysis and improved rates," {\it International Conference on Machine Learning}, 5200-5209,  2019.
 
    \bibitem{GraBoy:13}
	M. Grant and S. Boyd, "CVX: Matlab software for disciplined convex programming, version 2.0 beta," {\it http://cvxr.com/cvx}, 2013.

    \bibitem{LobLeb:98}
    M.S. Lobo, L. Vandenberghe, S. Boyd and H. Lebret, "Applications of second-order cone programming,"  {\it Linear Algebra and its Applications}, vol. 284, 193-228, 1998.
    \bibitem{LewPan:98}
	A. Lewis and J.S. Pang, \emph{Error bounds for convex inequality systems}, Generalized Convexity, Generalized Monotonicity (J.-P. Crouzeix, J.-E.Martinez-Legaz, and M. Volle, eds.), 75--110, Cambridge University Press, 1998.

    \bibitem{Ned:11}
	A. Nedich, "Random algorithms for convex minimization problems," {\it Mathematical Programming}, vol. 129, no. 2, 225--273, 2011.

    \bibitem{Nec:20}
	I. Necoara, "General convergence analysis of stochastic first order methods for composite optimization,"  {\it Journal of Optimization Theory and Applications}, vol. 189, 66-95,  2021.

    \bibitem{NecSin:22}
    I. Necoara and N.K. Singh, "Stochastic subgradient for composite  convex optimization with functional constraints," {\it Journal of Machine Learning Research}, vol. 23, no. 265, 1–35, 2022.


    \bibitem{NecNed:21}
    I. Necoara and A. Nedich, "Minibatch stochastic subgradient‑based projection
    algorithms for feasibility problems with convex inequalities," {\it Computational Optimization and Applications}, vol. 80, 121-152, 2021.

    \bibitem{NedDin:14}
    V. Nedelcu, I. Necoara and Q. Tran Dinh, "Computational complexity of inexact gradient augmented Lagrangian methods: application to constrained MPC", {\it SIAM Journal of  Control and Optimization}, vol. 52, no. 5, 3109–3134, 2014.

    \bibitem{Nes:18}
	Yu. Nesterov, "Lectures on Convex Optimization,"  {\it Springer Optimization and Its Applications}, vol. 137, 2018.

    \bibitem{Pol:69}
    B.T.  Polyak,  "Minimization of unsmooth functionals," {\it USSR Computational Mathematics and Mathematical Physics}, vol. 9, no. 3, 14-29, 1969.

    \bibitem{Pol:01}
    B.T.  Polyak, "Random algorithms for solving convex inequalities," {\it Studies in  Computational Mathematics}, vol. 8, 409-422,  2001.

    \bibitem{RocUry:00}
	R.T.  Rockafellar and S.P. Uryasev, "Optimization of conditional value-at-risk," {\it Journal of Risk},  vol. 2, 21-41, 2000.



    \bibitem{SidLuo:06}
    N.D. Sidiropoulos, T.N. Davidson and Z.-Q. Luo, "Transmit beamforming for physical-layer multicasting," {\it IEEE Transactions on Signal Processing}, vol. 54, no. 6, 2239–2251, 2006.

    \bibitem{SinNec:23}
    N.K. Singh, I. Necoara and V. Kungurtsev, "Mini-batch stochastic subgradient for functional constrained optimization," {\it Optimization}, 1–27, 2023, doi: 10.1080/02331934.2023.2189015.

    \bibitem{Tib:11}
	R. Tibshirani,  "The solution path of the generalized
		lasso," {\it Phd Thesis}, Stanford University, 2011.

    \bibitem{Vap:98}
	V. Vapnik, "Statistical learning theory," {\it John Wiley}, 1998.

    \bibitem{WanBer:16}
    M. Wang and D.P. Bertsekas, "Stochastic first-order methods with random constraint
    projection," {\it SIAM Journal on Optimization}, vol. 26, no. 1, 681–717, 2016.

    \bibitem{Wei:20}
    W. Wei, "Tutorials on advanced optimization methods," {\it 	arXiv:2007.13545}, 2020.

    \bibitem{Xu:20}
	Y. Xu, "Primal-dual stochastic gradient method for convex programs with many functional constraints," {\it SIAM Journal on Optimization}, vol. 30, no. 2, 1664–1692, 2020.

    \bibitem{YanYue:22}
	Z. Yang, F. Xia, K. Tu and M. Yue, "Variance reduced random relaxed projection method for constrained finite-sum minimization problems",  {\it IEEE Transactions on Signal Processing}, vol. 54, 2188-2203, 2024.
\end{thebibliography}
\end{document}